\documentclass{amsart}

\usepackage{amsfonts,amssymb,amscd,amsopn}
\usepackage{psfrag} \usepackage[dvips]{graphicx}
\usepackage[ansinew]{inputenc}
\usepackage{latexsym,psfrag,varioref,euscript,a4wide}

\usepackage{hyperref}  
  \hypersetup{colorlinks,
   debug=false,
    linkcolor=blue,  
    citecolor=red,  
    urlcolor=blue,   
   bookmarksopen=false,
   pdftitle={Chaotic vs singular-hyperbolic},
   pdfauthor={},
   pdfsubject={LaTeX},
   pdfkeywords={LaTeX}
}

\theoremstyle{plain}

\newtheorem{theorem}{Theorem }[section]
\newtheorem{proposition}[theorem]{Proposition}
\newtheorem{lemma}[theorem]{Lemma}
\newtheorem{corollary}[theorem]{Corollary}

\newtheorem{conjecture}{Conjecture}

\theoremstyle{definition} \theoremstyle{remark}
\newtheorem{remark}[theorem]{Remark}

\newtheorem{definition}[theorem]{Definition}

\newcommand{\ov}{\overline}

\newcommand{\leb}{\operatorname{Leb}}
\newcommand{\htop}{\operatorname{h_\textrm{top}}}

\newcommand{\dist}{\operatorname{dist}}

\newcommand{\inter}{\operatorname{int}}

\newcommand{\supp}{\operatorname{supp}}

\newcommand{\be} {\beta}        
       
\newcommand{\de} {\delta}       \newcommand{\De}{\Delta}

\newcommand{\la} {\lambda}      

\newcommand{\vfi}{\varphi}

\newcommand{\sgn}{\operatorname{sgn}}

\newcommand{\esssup}{\operatorname{ess\,sup}}

\def \CC {{\mathbb C}}
\def \DD {{\mathbb D}}
\def \EE {{\mathbb E}}

\def \II {{\mathbb I}}

\def \NN {{\mathbb N}}

\def \QQ {{\mathbb Q}}
\def \RR {{\mathbb R}}

\def \TT {{\mathbb T}}

\def \ZZ {{\mathbb Z}}

\newcommand{\cF}{\EuScript{F}}

\newcommand{\F}{\EuScript{F}}
\newcommand{\cD}{\EuScript{D}}
\newcommand{\cP}{\EuScript{P}}
\newcommand{\cO}{\EuScript{O}}

\newcommand{\U}{\EuScript{U}}

\newcommand{\X}{\EuScript{X}}
\newcommand{\cC}{\EuScript{C}}
\newcommand{\cS}{\EuScript{S}}

\newcommand{\B}{\EuScript{B}}
\newcommand{\R}{\EuScript{R}}

\newcommand{\varsq}{\operatorname{var^\square}}
\newcommand{\qand}{\quad\text{and}\quad}

\def \cA {{\mathcal A}}

\def \cC {{\mathcal C}}
\def \cD {{\mathcal D}}

\def \cF {{\mathcal F}}

\def \cH {{\mathcal H}}

\def \cO {{\mathcal O}}
\def \cP {{\mathcal P}}

\def \cR {{\mathcal R}}
\def \cS {{\mathcal S}}

\def \cU {{\mathcal U}}
\def \cV {{\mathcal V}}

\def \fX {{\mathfrak X}}

\def \fS {{\mathfrak S}}

\begin{document}

\title[Lorenz like flows, recent developments and
perspective] {Statistical properties of Lorenz like
  flows, recent developments and perspectives}

\thanks{This research has been supported in part by EU Marie-Curie IRSES
Brazilian-European partnership in Dynamical Systems (FP7-PEOPLE-2012-IRSES
318999 BREUDS).
V.A. and M.J.P. were partially supported by CNPq,
  PRONEX-Dyn.Syst., FAPERJ and FAPESB. M.J.P. was also
  partially supported by Balzan Research Project of J.Palis}

\author{Vitor Araujo}
\address[V.A.]{Instituto de Matem\'atica,
Universidade Federal da Bahia\\
Av. Adhemar de Barros, S/N , Ondina,
40170-110 - Salvador-BA-Brazil}
\email{vitor.d.araujo@ufba.br and vitor.araujo.im.ufba@gmail.com}

\author{Stefano Galatolo}
\address[S.G.]{Dipartimento di Matematica Applicata, Via Buonarroti 1 Pisa}
\email[S. Galatolo]{galatolo@dm.unipi.it}
\urladdr{http://users.dma.unipi.it/galatolo/}

\author{Maria Jos\'e Pacifico}
\address[M.J.P.]{Instituto de Matem\'atica,
Universidade Federal do Rio de Janeiro,\\
C. P. 68.530, 21.945-970 Rio de Janeiro, Brazil
}
\email{pacifico@im.ufrj.br and pacifico@impa.br}

\begin{abstract}
  We comment on mathematical results about the
  statistical behavior of Lorenz equations an its
  attractor, and more generally to the class of
  singular hyperbolic systems.  The mathematical theory
  of such kind of systems turned out to be surprisingly
  difficult. It is remarkable that a rigorous proof of
  the existence of the Lorenz attractor was presented
  only around the year 2000 with a computer assisted
  proof together with an extension of the hyperbolic
  theory developed to encompass attractors robustly
  containing equilibria.

  We present some of the main results on the
  statisitcal behavior of such systems.  We show that
  for attractors of three-dimensional flows, robust
  chaotic behavior is equivalent to the existence of
  certain hyperbolic structures, known as
  singular-hyperbolicity. These structures, in turn,
  are associated to the existence of physical measures:
  \emph{in low dimensions, robust chaotic behavior for
    flows ensures the existence of a physical measure}.

  We then give more details on recent results on the
  dynamics of singular-hyperbolic (Lorenz-like)
  attractors: (1) there exists an invariant foliation
  whose leaves are forward contracted by the flow (and
  further properties which are useful to understand the
  satistical properties of the dynamics); (2) there
  exists a positive Lyapunov exponent at every orbit;
  (3) there is a unique physical measure whose support
  is the whole attractor and which is the equilibrium
  state with respect to the center-unstable Jacobian;
  (4) this measure is exact dimensional; (5) the
  induced measure on a suitable family of
  cross-sections has exponential decay of correlations
  for Lipschitz observables with respect to a suitable
  Poincar\'e return time map; (6) the hitting time
  associated to Lorenz-like attractors satisfy a
  logarithm law; (7) the geometric Lorenz flow
  satisfies the Almost Sure Invariance Principle (ASIP)
  and the Central Limit Theorem (CLT); (8) the rate of
  decay of large deviations for the volume measure on
  the ergodic basin of a geometric Lorenz attractor is
  exponential; (9) a class of geometric Lorenz flows
  exhibits robust exponential decay of correlations;
  (10) all geometric Lorenz flows are rapidly mixing
  and their time-$1$ map satisfies both ASIP and CLT.
\end{abstract}

\keywords{sensitive dependence on initial conditions,
  physical measure, singular-hyperbolicity, expansiveness,
  robust attractor, robust chaotic flow, positive Lyapunov
  exponent, decay of correlations, large deviations,
  rapid mixing}

\date{\today}

\maketitle

\tableofcontents

\section{Introduction}
\label{sec:introd}

\begin{figure}[htpb]
\psfrag{X}{$X^t(x)$}\psfrag{Y}{$X^t(y)$}
\psfrag{x}{$x$}\psfrag{y}{$y$}
\includegraphics[width=4cm]{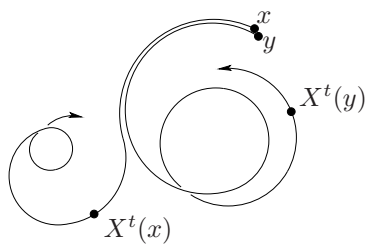}
\caption{\label{fig-sensivel} Sensitive dependence on
  initial conditions.}
\end{figure}

The development of the theory of dynamical systems has
shown that many important dynamical models exhibit
\emph{sensitive dependence on initial
  conditions}\footnote{Formally the definition of
  sensitivity for a flow $X^t$ on some compact manifold
  $M$ is as follows: an $X^t$-invariant subset $\Lambda$
  is \emph{sensitive to initial conditions} or has
  \emph{sensitive dependence on initial conditions}, or
  simply \emph{chaotic} if, for every small enough
  $r>0$ and $x\in\Lambda$, and for any neighborhood $U$
  of $x$, there exists $y\in U$ and $t\neq0$ such that
  $X^t(y)$ and $X^t(x)$ are $r$-apart from each other:
  $\dist\big(X^t(y),X^t(x)\big)\ge r$; see
  Figure~\ref{fig-sensivel} and
  Section~\ref{sec:robust-chaotic-parti}. An analogous
  definition holds for maps $f$ of some manifold or
  even metric spaces.
}, a common feature of \emph{chaotic dynamics}: small
initial differences are rapidly augmented as time
passes, causing two trajectories originally coming from
practically indistinguishable points to behave in a
completely different manner after a short
while. Pointwise, long term predictions based on such
models are unfeasible since it is not possible to both
specify initial conditions with arbitrary accuracy and
numerically calculate with arbitrary precision; for an
introduction to these notions
see~\cite{devaney1989,robinson2004}.

On the other hand in these systems, even if the pointwise description or forecasting of the system is forbidden by the initial condition sensitivity, the statistical behavior is often relatively simple and its properties are often (with a certain effort) predictable.

A theory which explain this statistical behavior is quite satisfactorily developped for systems having some uniformly hyperbolic behavior (see below for the definition), yet the rigorous description of the statistical behavior of relatively simple systems as: quadratic polynomials \footnote{(as the \emph{logistic family} or \emph{H\'enon
  attractor}, see e.g.\cite{devaney1989} for a gentle
introduction)}, or autonomous ordinary differential equations with a hyperbolic equilibrium of saddle-type accumulated by regular orbits, as the \emph{Lorenz flow}  is still far from being complete. In this article we are going to describe some relatively recent developments in this second direction.

\subsection{The Lorenz equations}
\label{sec:lorenz-equati}

In 1963 the meteorologist Edward Lorenz published in the
Journal of Atmospheric Sciences \cite{Lo63} an example of a
parametrized polynomial system of differential equations
\begin{align}
  \label{e-Lorenz-system}
\dot x &= a(y - x)&
\quad
&a = 10 \nonumber
\\
\dot y &= rx -y -xz&
\quad
&r =28
\\
\dot z &= xy - bz&
\quad
&  b = 8/3 \nonumber
\end{align}
as a very simplified model for thermal fluid convection,
motivated by an attempt to understand the foundations of
weather forecast.

The origin $\sigma=(0,0,0)$ is an equilibrium of saddle type for the vector 
field defined by equations (\ref{e-Lorenz-system}) with real eigenvalues 
$\lambda_i$, $i\leq 3$ satisfying
\begin{equation}
\label{eigenvalues}
  \lambda_2<\lambda_3<0<-\lambda_3<\lambda_1.
\end{equation}
(in this case $\lambda_1\approx 11.83$ , $\lambda_2\approx
-22.83$, $\lambda_3=-8/3$).

\begin{figure}[htpb]
  \includegraphics[scale=0.7]{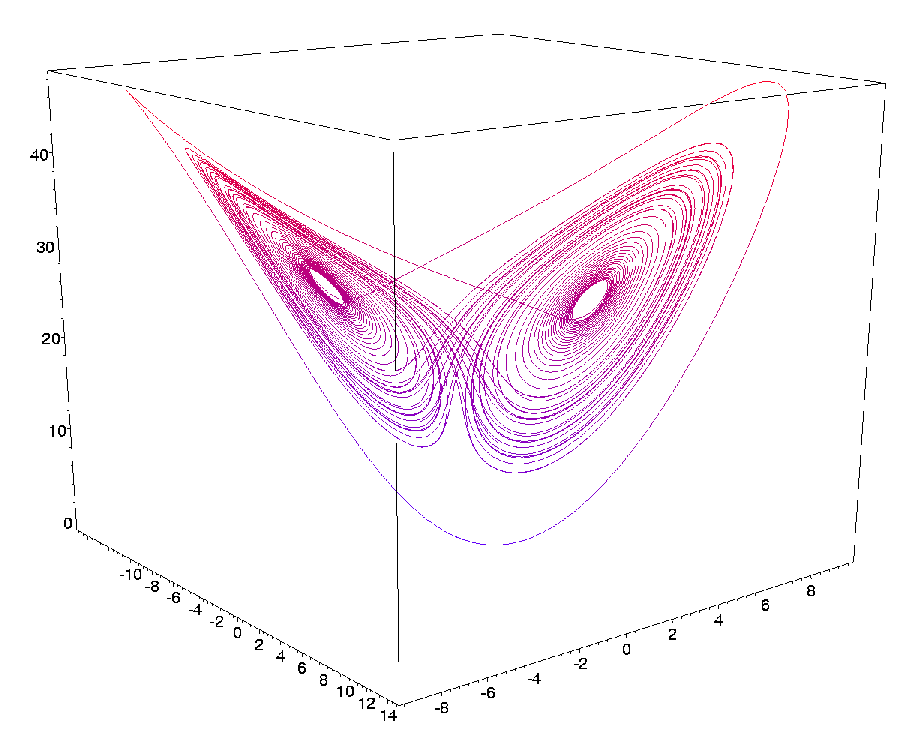}
  \caption{A view of the Lorenz attractor calculated numerically}
  \label{fig:view-lorenz-attract}
\end{figure}

Numerical simulations performed by Lorenz for an open
neighborhood of the chosen parameters suggested that almost
all points in phase space tend to a {\em chaotic attractor},
that is, a bounded region in phase-space, invariant under
time evolution, such that the forward trajectories of most
or even all points nearby converge to it, and these
trajectories are {\em sensitive with respect to initial data}.
The well known
picture of the Lorenz attractor is presented in
Figure~\ref{fig:view-lorenz-attract}.

Lorenz's equations proved to be very resistant to
rigorous mathematical analysis, and also presented
serious difficulties to rigorous numerical study.
Indeed, these two main difficulties are:
\begin{description}
\item[\emph{conceptual}] the presence of an equilibrium point
  at the origin accumulated by regular orbits of the flow
  prevents this attractor from being hyperbolic
  \cite{AraPac2010},
\item[\emph{numerical}] the presence of an equilibrium
  point at the origin, implying that solutions slow down as
  they pass near the origin, which means unbounded return
  times and, thus, unbounded integration errors.
\end{description}
Moreover the attractor is \emph{robust}, that is, the
features of the limit set persist for all nearby vector
fields. More precisely, if $U$ is an isolating neighborhood
of the attractor $\Lambda$ for a vector field $X$, then
$\Lambda$ is \emph{robustly transitive} if, for all vector
fields $Y$ which are $C^1$ close to $X$, the corresponding
$Y$-invariant set
\begin{align*}
  \Lambda_Y(U)=\bigcap_{t>0} Y^t(U)
\end{align*}
also admits a dense positive $Y$-orbit.  The
\emph{persistence of transitivity}, that is, the fact that,
for all nearby vector fields, the corresponding limit set is
transitive, \emph{is quite remarkable and implies a
  dynamical characterization of the attractor}, as we shall
see.

\subsection{Geometric Lorenz model, computer aided approach, singular hyperbolic attractors}
\label{sec:geometr-lorenz-model}

The difficulties in the study of the Lorenz system led, in
the seventies, to the construction of geometric flows
presenting a similar behavior as the one generated by
equations (\ref{e-Lorenz-system}). Nowadays these models are
known as {\em geometric Lorenz flows}.  We describe this
construction in Section~\ref{sec:geoLorenz}; see
\cite{ACS95, GW79} for full details.

These models are three-dimensional flows for which it is
easy to rigorously prove the existence of an attractor
containing an equilibrium point of the flow, together with
regular solutions. 
They have a natural  cross-section given by a
two-dimensional square crossed by all orbits of the flow
inside the attractor except the singularity. 
Several properties about the statistical behavior of these flows have been rigorously understood by the use of  the properties   of the   Poincar\'e first return map to the cross-section.

A succesful approach for the real Lorenz flow was through
rigorous numerics. In this way, it could be proved
\cite{HHTZ94,HT92,MM95,MM98} that the Lorenz system of
equations exhibits a suspended Smale horseshoe \cite{Sm67}
which implies, in particular, the existence of infinitely
many closed orbits. However, proving the existence of an
attractor as in the geometric models is an even harder task,
because one cannot avoid the fact that solutions slow down
as they pass near the equilibrium, which, as said before
means unbounded return times and so unbounded integration
errors.  This was finally settled by Tucker in
\cite{Tu2} around the turn of the century.

Using a combination of rigorous numerics and normal form
theory, Tucker proved that the Lorenz equations
\eqref{e-Lorenz-system} support a robust strange attractor
$\Lambda $ and the flow admits a unique Sinai-Ruelle-Bowen measure
$\mu$ with $\supp(\mu)=\Lambda $; this notion will be presented in
Section~\ref{sec:prelim-notions}.  Tucker's proof uses a
computer algorithm to estimate convenient solutions of
\eqref{e-Lorenz-system}, keeping rigorous bounds on the
errors. Successful termination of this algorithm proves the
presence of a robustly transitive attractor in
\eqref{e-Lorenz-system}.
It is worth to remark that even in this approach a key step is the construction of a suitable Poincar\'e map on a cross section and the understanding of its properties (see Section \ref{3} ) which are similar to the ones found in the geometric model.

From robust transitivity it follows after \cite{MPP98,MPP04}
that the attractor supported by the equations
\eqref{e-Lorenz-system}  is a singular-hyperbolic attractor.
We will see in the following that also for these attractors we can construct a suitable cross section and prove that they shares all the fundamental features of the
geometric Lorenz models (see Section \ref{sec:singul-hyperb-system} )  so that its geometry and its
ergodic properties can be well understood.

Singular hyperbolicity plays the role of hyperbolicity for
flows presenting equilibria accumulated by regular
orbits. Hyperbolicity alone, a classical notion going back
to Smale \cite{Sm67} means that the complementary direction
to the flow can be further split into a
pair of complementary invariant directions, one uniformly
contracting and the other uniformly expanding
by the tangent map to the flow.  The theory
of hyperbolic systems describing their geometric and ergodic
properties is very rich: the interested reader should
consult \cite{BR75,GW79,PM82,Sh87} and references therein.

Singular hyperbolicity replaces the expanding direction by a
two-dimensional direction containing the flow direction on
regular orbits along which the flow should expand area.
Singular hyperbolicity encompasses flows exhibiting
equilibria attached to regular orbits, since by \cite{MPP04}
an attractor robustly containing equilibria is singular
hyperbolic.

Moreover singular hyperbolicity is a natural generalization
of hyperbolicity, since compact invariant sets which are
singular hyperbolic, but have no equilibria, can be proved
to be hyperbolic in the usual sense.  Remarkable
dynamical properties can be proved for singular hyperbolic attractors
extending in this way the hyperbolic theory to a wider class
of systems. We present some of these dynamical and ergodic
properties of singular-hyperbolic attractors in what
follows.

\subsection{Overview}
\label{sec:overvi}

After reviewing some main ideas from dynamical systems
theory, using results on robustness of attractors from
Ma\~n\'e~\cite{Man82} and Morales, Pacifico and
Pujals~\cite{MPP04} together with observations on their
proofs, we show that for attractors of three-dimensional
flows, robust chaotic behavior (in the above sense of
sensitiveness to initial conditions for all
  close enough flows) is equivalent to the existence of
certain partially hyperbolic structures. These structures,
in turn, allow the construction of suitable Poincar\'e sections with nice properties and to deduce several consequences about the statistical behavior of the dynamics.

In the following, we review the construction and
several more or less recent results about the dynamics
of singular-hyperbolic (or Lorenz-like) attractors:
\begin{itemize}
\item there exists an invariant foliation whose leaves
  are forward contracted by the flow and the dynamics
  satisfies a list of properties which is similar to
  the ones of geometric Lorenz attractors; see
  Section~\ref{sec:robustn-singul-hyper}.
\item there exists a positive Lyapunov exponent at every
  orbit; see Section~\ref{sec:robustn-singul-hyper}.
\item there is a unique physical measure whose support is
  the whole attractor and which is the equilibrium state
  with respect to the center-unstable Jacobian; see
  Section~\ref{sec:ergodic-theory-singu}.
\item this physical measure is exact dimensional and
  the hitting time associated to a Lorenz-like
  attractor satisfies a logarithm law; see
  Sections~\ref{sec:muexata} and
  \ref{sec:logarithm-law-syngul}.
\item the induced measure on a suitable family of
  cross-sections of a Lorenz-like flow has exponential decay
  of correlations with respect to the Poincar\'e first
  return map; see Section~\ref{sec:decay-correl-poincar-1}.
\item the geometric Lorenz flow satisfies the Almost Sure
  Invariance Principle and the Central Limit Theorem; see
  Section~\ref{sec:central-limit-theore}.
\item the rate of decay for large deviations with respect to
  the volume measure on the ergodic basin of a geometric
  Lorenz attractor is exponential; see Section~\ref{sec:large-deviations}.
\item there are open sets of geometric Lorenz flows
  each of which exhibits exponential decay of
  correlations; see
  Section~\ref{sec:exponent-decay-corre-geomLorenzflow}.

\item In low dimensions, robust chaotic behavior
  ensures the existence of a physical measure; see
  Section~\ref{sec:chaotic-systems}. \footnote{We
    remark that this provides a partial answer to a
    conjecture of Viana: \emph{existence of positive
      Lyapunov exponent for a positive Lebesgue measure
      subset of orbits implies the existence of a
      physical measure} This is proved in
    Section~\ref{sec:robust-chaotic-parti}}

\end{itemize}

We finish with a brief list of conjectures on the dynamics
of singular-hyperbolic attractors, in
Section~\ref{sec:some-conject-open}.


\section{Preliminary notions}
\label{sec:prelim-notions}

Here and throughout the text we assume that $M$ is a
three-dimensional compact connected manifold without
boundary endowed with some Riemannian metric which induces a
distance denoted by $\dist$ and a volume form $\leb$ which
we name \emph{Lebesgue measure} or \emph{volume}. For any
subset $A$ of $M$ we denote by $\ov{A}$ the (topological)
closure of $A$.

We denote by $\fX^r(M), r\ge1$ the set of $C^r$ smooth
vector fields $X$ on $M$ endowed with the $C^r$
topology. Given $X\in\fX^r(M)$ we denote by $X^t$, with
$t\in\RR$, the flow generated by the vector field $X$. Since
we assume that $M$ is a compact manifold the flow is defined
for all time. Recall that the flow
$(X^t)_{t\in\RR}$ is a family of $C^r$ diffeomorphisms
satisfying the following properties:
\begin{enumerate}
\item $X^0=Id:M\to M$ is the identity map of $M$;
\item $X^{t+s}=X^t\circ X^s$ for all $t,s\in\RR$,
\end{enumerate}
and it is \emph{generated by the vector field} $X$ if
\begin{enumerate}
\item[(3)] $\left.\frac{d}{dt} X^t (q)\right|_{t=t_0} =
X\big(X_{t_0}(q)\big)$ for all $q\in M$ and $t_0\in\RR$.
\end{enumerate}

We say that a compact $X^t$-invariant set $\Lambda$ is
\emph{isolated} if there exists a neighborhood $U$ of
$\Lambda$ such that $\Lambda=\cap_{t\in\RR}X^t(U)$.  A
compact invariant set $\Lambda$ is {\em attracting} if
$\Lambda_X(U):=\cap_{t\geq 0}X^t(U)$ equals $\Lambda$ for
some neighborhood $U$ of $\Lambda$ satisfying
$\ov{X^t(U)}\subset U$, for all $t> 0$. In this case the
neighborhood $U$ is called an \emph{isolating} neighborhood
of $\Lambda$. Note that $\Lambda_X(U)$ is in general
different from $\cap_{t\in\RR} X^t(U)$, but for an
attracting set the extra condition $\ov{X^t(U)}\subset U$
for $t>0$ ensures that every attracting set is also
isolated. We say that $\Lambda$ is {\em transitive} if
$\Lambda$ is the closure of both $\{X^t(q):t>0\}$ and
$\{X^t(q):t<0\}$ for some $q\in \Lambda$.  An {\em
  attractor} of $X$ is a transitive attracting set of $X$
and a {\em repeller} is an attractor for $-X$.  We say that
$\Lambda$ is a \emph{proper} attractor or repeller if
$\emptyset\neq\Lambda\neq M$.

An \emph{equilibrium} (or \emph{singularity}) for $X$ is a
point $\sigma\in M$ such that $X^t(\sigma)=\sigma$ for all
$t\in\RR$, i.e.  a fixed point of all the flow maps, which
corresponds to a zero of the associated vector field $X$:
$X(\sigma)=0$. An \emph{orbit} of $X$ is a set
$\cO(q)=\cO_X(q)=\{X^t(q): t\in\RR\}$ for some $q\in M$.  A
{\em periodic orbit} of $X$ is an orbit $\cO=\cO_X(p)$ such
that $X^T(p)=p$ for some minimal $T>0$. A \emph{critical
  element} of a given vector field $X$ is either an
equilibrium or a periodic orbit.

We say that a compact invariant subset is \emph{singular
  hyperbolic} if all the singularities in $\Lambda$ are
hyperbolic, and the tangent bundle $T\Lambda$ decomposes in
two complementary $DX^t$-invariant bundles $E^s \oplus
E^{cu}$, where: $E^s$ is one-dimensional and uniformly
contracted by $DX^t$; $E^{cu}$ is bidimensional, contains the
flow direction, $DX^t$ expands area along $E^{cu}$ and
$DX^t\mid E^{cu}$ dominates $DX^t\mid E^s$ (i.e.  any eventual
contraction in $E^s$ is stronger than any possible
contraction in $E^{cu}$), for all $t>0$.  

The notion of singular hyperbolicity was introduced
in~\cite{MPP98,MPP04} where it was proved that any $C^1$ robustly
transitive set for a $3$-flow is either a singular
hyperbolic attractor or repeller.

We note that the presence of an equilibrium together with
regular orbits accumulating on it prevents any invariant set
from being  hyperbolic, see
e.g.~\cite{BR75}. Indeed, in our $3$-dimensional setting a
compact invariant subset $\Lambda$ is
 hyperbolic if the tangent
bundle $T\Lambda$ decomposes in \emph{three} complementary
$DX^t$-invariant bundles $E^s \oplus E^{X} \oplus E^u$, each
one-dimensional, $E^X$ is the flow direction, $E^s$ is
uniformly contracted and $E^u$ uniformly expanded by $DX^t$,
$t>0$. This implies the continuity of the splitting and the
presence of a non-isolated equilibrium point in $\Lambda$
leads to a discontinuity in the splitting dimensions.

In the study of the asymptotic behavior of orbits of a flow
$X \in {\fX}^1(M)$, a fundamental problem is to understand
how the behavior of the tangent map $DX$ determines the
dynamics of the flow $X^t$.  The main achievement along this
line is the uniform hyperbolic theory: we have a very good
description of the dynamics assuming that the tangent map
has a hyperbolic structure since the work of Bowen and
Ruelle~\cite{BR75}.

We recall standard facts about hyperbolic flows from
e.g. \cite{HPS77}. An embedded disk $\gamma\subset M$ is a
(local) {\em strong-unstable manifold}, or a {\em
  strong-unstable disk}, if $\dist(X^{-t}(x),X^{-t}(y))$
tends to zero exponentially fast as $t\to+\infty$, for every
$x,y\in\gamma$. Similarly, $\gamma$ is called a (local) {\em
  strong-stable manifold}, or a {\em strong-stable disk}, if
$\dist(X^{t}(x),X^{t}(y))\to0$ exponentially fast as
$n\to+\infty$, for every $x,y\in\gamma$. It is well-known
that every point in a hyperbolic set possesses a local
strong-stable manifold $W_{loc}^{ss}(x)$ and a local
strong-unstable manifold $W_{loc}^{uu}(x)$ which are disks
tangent to $E_x$ and $G_x$ at $x$ with topological
dimensions $d_s=\dim(E^s)$ and $d_u=\dim(E^u)$,
respectively. 
These disks are $X^t$-invariant, meaning for $x\in\Lambda$
and $t>0$
\begin{align*}
  X^t(W^{ss}_{loc}(x))\subset W^{ss}_{loc}(X^t(x))
  \qand
  X^{-t}(W^{uu}_{loc}(x))\subset W^{uu}_{loc}(X^{-t}(x))
\end{align*}
and so we obtain the (global) \emph{strong-stable manifold}
$$W^{ss}(x)=\bigcup_{t>0}
X^{-t}\Big(W^{ss}_{loc}\big(X^t(x)\big)\Big)$$
and the
(global) \emph{strong-unstable manifold}
$$W^{uu}(x)=\bigcup_{t>0}X^{t}\Big(W^{uu}_{loc}\big(X^{-t}(x)\big)\Big)$$
for every point $x$ of a hyperbolic set.  These are immersed
submanifolds with the same differentiability of the flow.
We also consider the \emph{stable manifold}
$W^s(x)=\cup_{t\in\RR} X^{t}\big(W^{ss}(x)\big)$ and
\emph{unstable manifold}
$W^u(x)=\cup_{t\in\RR}X^{t}\big(W^{uu}(x)\big)$ for $x$ in a
hyperbolic set, which are flow invariant.

In the same vein, under the assumption of singular
hyperbolicity and also following the standard reference
\cite{HPS77}, one can show that at each point there exists a
strong stable manifold and that the whole set is foliated by
leaves that are contracted by forward iteration.

In particular, this shows that any robust transitive
attractor with singularities displays similar properties to
those of the geometric Lorenz model.  It is also possible
to show the existence of local central manifolds tangent to
the central unstable direction; see e.g. \cite{ArGalPac}.
Although these central manifolds do not behave as unstable
ones, in the sense that points on them are not necessarily
asymptotic in the past. The expansion of volume along the
central unstable two-dimensional direction enables us to
deduce some remarkable properties.

We recall that a $X^t$-invariant probability measure $\mu$
is a probability measure satisfying $\mu(X^t(A))=\mu(A)$ for
all $t\in\RR$ and measurable $A\subset M$. Given an
invariant probability measure $\mu$ for a flow $X^t$, let
$B(\mu)$ be  the \emph{(ergodic) basin} of $\mu$, i.e.,
the set of points $z\in M$ satisfying for all continuous
functions $\varphi: M \to \RR$
\begin{align*}
  \lim_{T\to+\infty} \frac{1}{T} \int_0^T \varphi\big(
  X^t(z) \big) \, dt = \int\varphi\,d\mu.
\end{align*}
We say that $\mu$ is a \emph{physical} (or \emph{SRB})
measure for $X$ if $B(\mu)$ has positive Lebesgue measure:
$\leb\big( B(\mu) \big)>0$. 

The existence of a physical measures for an attractor shows
that most points in a neighborhood of the attractor have
well defined long term statistical behavior. So, in spite of
chaotic behavior preventing the exact prediction of the time
evolution of the system in practical terms, we gain some
statistical knowledge of the long term behavior of the
system near the chaotic attractor.

\section{Geometric Lorenz system(s)}
\label{3}

In this section we introduce a concrete flow which is in
some sense the simplest example of singular-hyperbolic
system. This model was also historically, the first one
where some rigorous results on the dynamic of Lorenz like
system were proved.  The system has a fixed point at the
origin, and a linear vector field around it. Its global
behavior is similar to the original Lorenz system but the
linearity near the origin allows to write explicit formulas
for the trajectory near it, and an explicit form for the
Poincar\'e map on a suitable section. This allow to obtain
many properties of the flow, which can be used to deduce
several statistical consequences for its dynamics.
We remark that,  sometimes in the literature for geometric Lorenz system is
meant a system satisfying a list of properties as the ones
we will see in the geometric Lorenz one:
\begin{enumerate}
\item there exists a Lorenz-like singularity at the origin for a
  $C^2$ smooth vector field in $\RR^3$;
\item there is a suitable cross-section $\Sigma$ whose
  Poincar\'e first return map $P$ preserves a uniformly
  contracting fibration;
\item the one dimensional induced map $f$ on the quotient
  $X$ by this contracting fibration is piecewise $C^{1+\epsilon}$, for
  some $\epsilon>0$, with two branches and has a singularity
  corresponding to a certain leaf $\xi_0$: if $d(x)$ denotes
  the distance of $x\in\Sigma$ to $\xi_0$, then there exists
  $\beta\in(0,1)$ such that $|f'(x)|=d(x)^{\beta-1}g(x)$
  with $g\in C^\epsilon(\Sigma)$;
\item the Poincar\'e first return time $r:\Sigma\to\RR^+$ is
  integrable with respect to Lebesgue area measure on the
  cross-section and there exists a constant $c_0>0$ such
  that $r(x)=-c_0\log d(x)+h(x)$ with $h\in
  C^\epsilon(\Sigma)$.
\item $f$ is \emph{uniformly piecewise expanding}: there are constants
  $\sigma>1$ and $c>0$ such that $|(f^n)'(x)|\ge c\sigma^n$
  for all $x\in X$ and $n>1$.
\end{enumerate}
The precise assumptions which are considered may slightly
vary from paper to paper. 
To avoid confusion, in this paper we will refer to a systems considered in this approach as
\emph{Axiomatic geometric Lorenz system.}  Systems of this
kind however do behave in a similar way having the same
asymptotic features.

\subsection{The construction}
\label{sec:geoLorenz}

Next we briefly recall the construction of a concrete example of  a flow satisfying the above list of properties.
We  will refer this as the geometric 
Lorenz flow, that is the simpler example of
a singular-hyperbolic attractor; see \cite{ACS95, GW79} for
full details.  As explained in the Introduction, the purpose
was the construction of a geometric flow presenting a
similar behavior as the one generated by equations
(\ref{e-Lorenz-system}). 
We start by observing that under some non-resonance
conditions, by the results of Sternberg~\cite{St58}, in a
neighborhood of the origin, which we assume to contain the
cube $[-1,1]^3 \subset \RR^3$, the Lorenz equations are
equivalent to the linear system $(\dot x, \dot y, \dot
z)=(\lambda_1 x,\lambda_2 y, \lambda_3 z)$ through
smooth conjugation, thus
\begin{align}\label{eq:LinearLorenz}
X^t(x_0,y_0,z_0)=
(x_0e^{\lambda_1t}, y_0e^{\lambda_2t}, z_0e^{\lambda_3t}),
\end{align}
where $\lambda_1\approx 11.83$ , $\lambda_2\approx -22.83$,
$\lambda_3=-8/3$ and $(x_0, y_0, z_0)\in\RR^3$ is an
arbitrary initial point near $(0,0,0)$.

Consider $S=\big\{ (x,y,1) : |x|\le
{\scriptstyle{1/2}},\quad |y|\le{\scriptstyle{1/2}}\big\}$ and
\begin{align*}
S^-&=\big\{ (x,y,1)\in S : x<0 \big\},&
\qquad
S^+&=\big\{ (x,y,1)\in S : x>0 \big\}\quad\text{and}
\\
S^*&=S^-\cup S^+=S\setminus\ell, & \text{where}\quad
\ell &=\big\{(x,y,1)\in S : x=0 \big\}.
\end{align*}
Assume that $S$ is a global transverse section to the flow so that
every trajectory eventually crosses $S$ in the direction of
the negative $z$ axis.

Consider also $\Sigma=\{ (x,y,z) : |x|=1
\}={\Sigma}^-\cup{\Sigma}^+$ with ${\Sigma}^{\pm}=\{ (x,y,z)
: x=\pm 1\}$. 

For each $(x_0,y_0,1)\in S^*$ the time $\tau$
such that $X^{\tau}(x_0,y_0,1)\in\Sigma$ is given by
$$\tau(x_0)=-\frac{1}{\lambda_1}\log{|x_0|},$$ which depends
on $x_0\in S^*$ only and is such that $\tau(x_0)\to+\infty$
when $x_0\to 0$. This is one of the reasons many standard
step by step numerical integration algorithms were unsuited to tackle the Lorenz
system of equations. Hence we get
(where $\sgn(x)=x/|x|$ for $x\neq0$)
\begin{equation}\label{L} 
X^\tau(x_0,y_0,1)=
\big( \sgn(x_0),  y_0e^{\lambda_2\tau},
e^{\lambda_3\tau}\big)
=
\big( \sgn(x_0),
y_0|x_0|^{-\frac{\lambda_2}{\lambda_1}},
|x_0|^{-\frac{\lambda_3}{\lambda_1}}\big).
\end{equation}
Since $0<-\lambda_3<\lambda_1<-\lambda_2$, we  have
$0<\alpha=-\frac{\lambda_3}{\lambda_1} <1
<\beta=-\frac{\lambda_2}{\lambda_1}$.
Let $L:S^*\to\Sigma$ be such that $ L(x,y)=\big(
y|x|^\beta,|x|^\alpha \big)$ with the convention that
$L(x,y)\in{\Sigma}^+$ if $x>0$ and $L(x,y)\in{\Sigma}^-$ if $x<0$.
\begin{figure}[h]
\includegraphics[width=7.5cm]{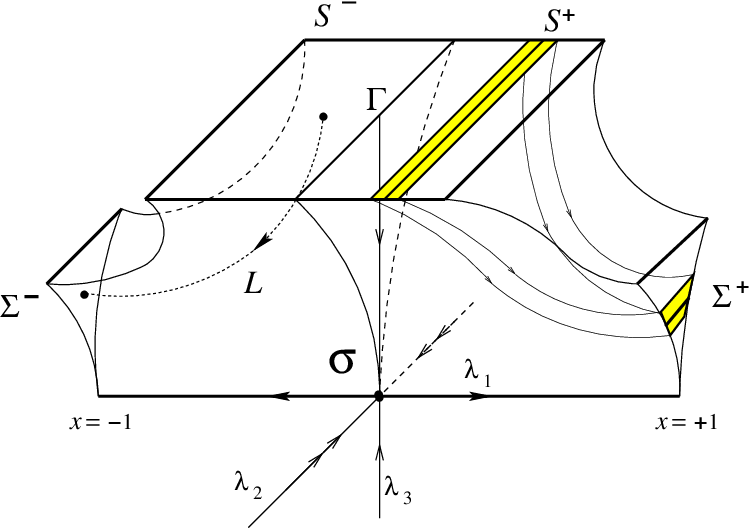}
\caption{\label{L3Dcusp}Behavior near the origin.}
\end{figure}
It is easy to see that $L(S^\pm)$ has the shape of a
 cusp like triangle without the vertex $(\pm 1,0,0)$.
In fact the vertex $(\pm 1,0,0)$ are
cusp points at the boundary of each of these sets.
The fact that $0<\alpha<1<\beta$ together with equation (\ref{L}) imply
that $L(\Sigma^\pm)$ are uniformly compressed in the $y$-direction.

From now on we denote by $\Sigma^\pm$ the closure of
$L(S^\pm)$.  Clearly each line segment $S^*\cap\{x=x_0\}$ is
taken to another line segment $\Sigma\cap\{z=z_0\}$ as
sketched in Figure~\ref{L3Dcusp}.

The sets $\Sigma^\pm$ should return to the cross section $S$
through a composition of a family of translations $T_t$, a
family of expansions $E_t$ only along the $x$-direction and
a family of rotations $R_t$ around $W^s(\sigma_1)$ and
$W^s(\sigma_2)$, where $\sigma_i$ are saddle-type
singularities of $X^t$ that are outside the cube $[-1,1]^3$,
see \cite{AraPac2010}. We assume that this composition takes
line segments $\Sigma\cap\{z=z_0\}$ into line segments
$S\cap\{x=x_1\}$ as sketched in Figure~\ref{L3Dcusp}.  The
composition $T_t \circ E_t \circ R_t $ of linear maps
describes a vector field $V$ in a region $W$ outside $[-1,1]^3$.
The geometric Lorenz flow $X^t$ is then defined in the
following way: for each $t\in\RR$ and each point $x \in S$,
the orbit $X^t(x)$ will start following the linear field
until $\tilde\Sigma^\pm$ and then it will follow $V$ coming
back to $S$ and so on. Let us write $\B=\{ X^t(x), x\in S,
t\in \RR^+\} $ the set where this flow acts.  The geometric
Lorenz flow is then the pair $({\B}, X^t )$ defined in this
way.  The set
$$
\Lambda=\cap_{t\ge 0}X^t(W\cup[0,1]^3)
$$
is the {\em geometric Lorenz attractor}.

\begin{figure}[htb]
  \centering
  \psfrag{e}{$\ell$}\psfrag{y}{$f(x)$}\psfrag{S}{$S$}
  \psfrag{w}{$-1$}\psfrag{z}{$0$}\psfrag{u}{$1$}
  \psfrag{l1}{$\lambda_1$}\psfrag{l2}{$\lambda_2$}\psfrag{l3}{$\lambda_3$}
  \includegraphics[width=10cm]{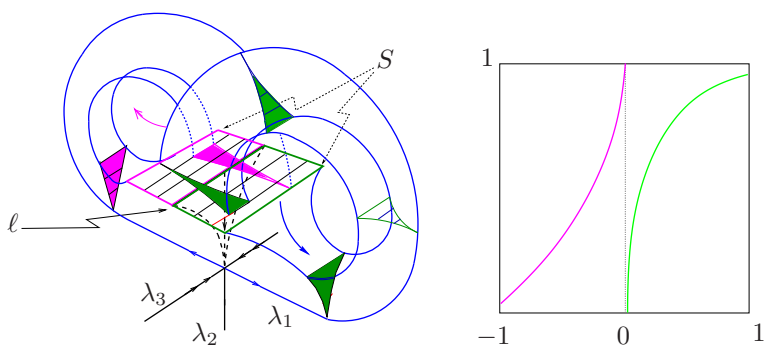}
\caption{\label{fig:GeoLorenz} The global cross-section for the geometric Lorenz
  flow and the associated 1d quotient map, the Lorenz transformation.}
\end{figure}

The combined effects of $T\circ E\circ R$ and the linear
flow given by equation (\ref{L}) on lines implies that the
foliation $\cF^s$ of $S$ given by the lines $S\cap\{x=x_0\}$
is invariant under the first return map $F:S\setminus\ell\to
S$. In other words, we have {\em for any given leaf $\gamma$
  of $\cF^s$, its image $F(\gamma)$ is contained in a leaf
  of $\cF^s$}.

\begin{figure}[htbp]
\includegraphics[width=5cm]{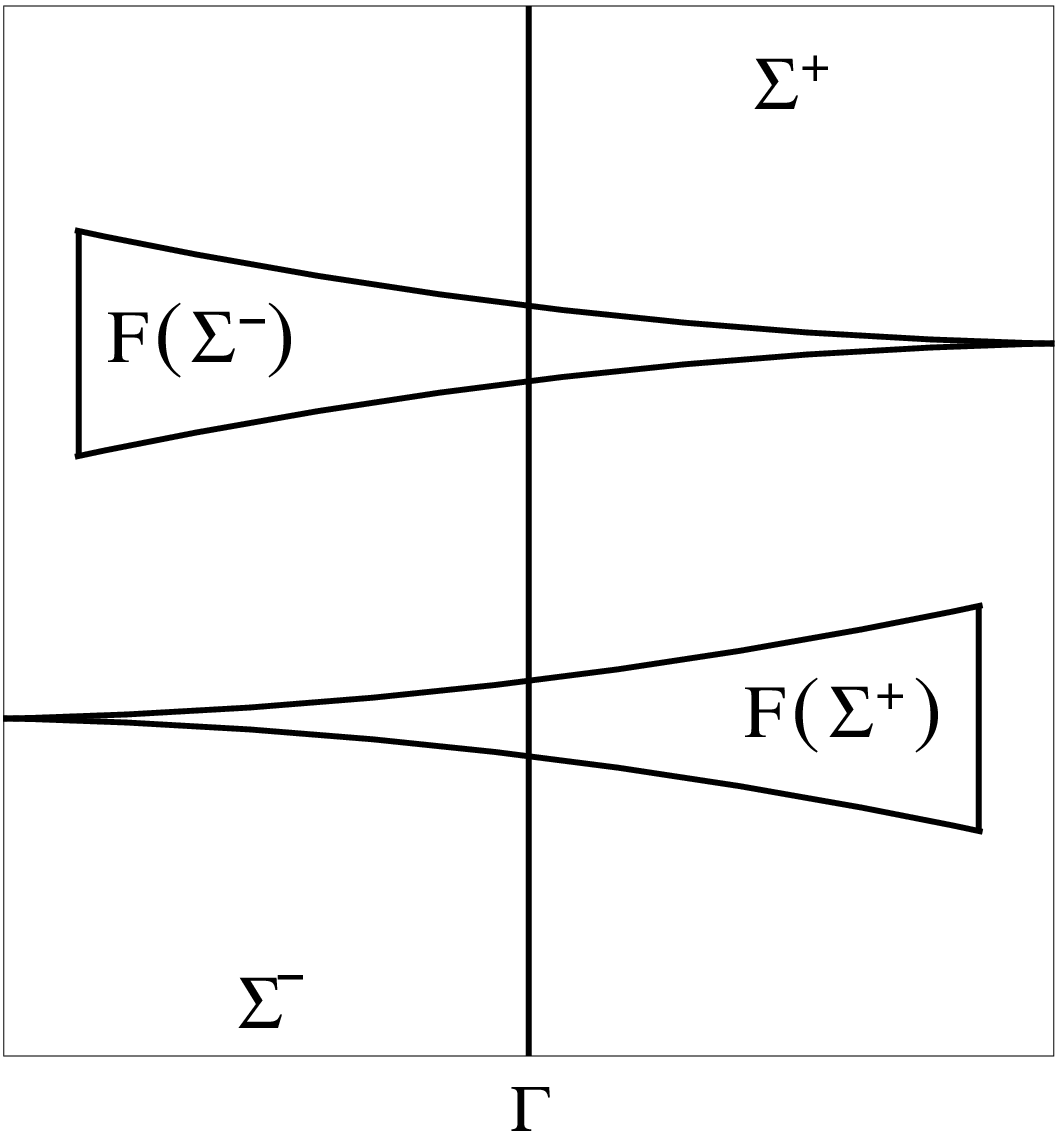}
\caption{\label{L2D}{The image $F(S^*)$.}}
\end{figure}

The main features of the geometric Lorenz flow and its first return map
can be seen at figures \ref{fig:GeoLorenz} and \ref{L2D}.

\begin{figure}[htbp]
\includegraphics[width=5cm]{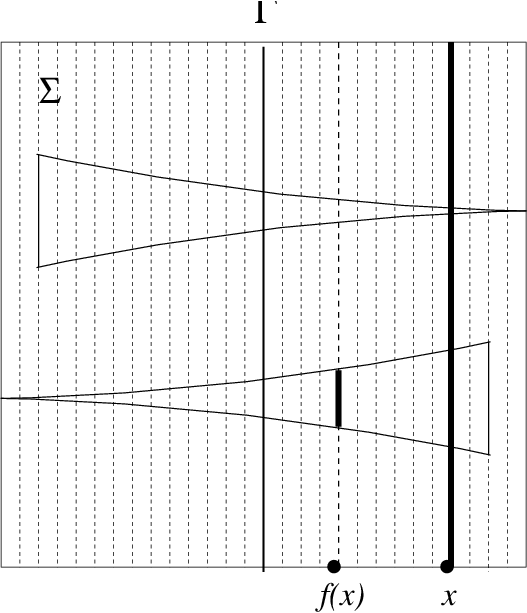}
\caption{\label{folhea}{Projection on $I$.}}
\end{figure}

The invariance of the foliation $\cF^s$ by lines
$S\cap\{x=x_0\}$ ensures that the Poincar\'e first return
map $F:S\setminus\ell\to S$ can be written as a skew-product
$F(x,y)=(f(x),g(x,y))$ and the one-dimensional map
$f:[-1/2,0)\cup(0,1/2]\to[1/2,1/2]$ is also the quotient map
of $F$ over the leaves of the stable foliation $\cF^s$
defined above. 
It is crucial that this map
  is \emph{piecewise expanding} and has a singularity at
  zero. More precisely, $f$ is $C^1$ on each
  interval $[-1/2,0), (0,1/2]$, its derivative is
  H\"older continuous and satisfies
  \begin{itemize}
  \item $f'(x)=|x|^\beta h(x)$ for
    $\alpha=-\lambda_3/\lambda_1$ and $h:[-1/2,1/2]\to\RR$
    H\"older continuous (note that $f'(0^+)=f'(0^-)=+\infty$).
  \item $f'$ is uniformly expanding, that is, there are
  $C>0$ and $\lambda>1$ such that $|(f^n)'|>C\lambda^n$ for
  all iterates $n\ge1$ and at all points where it is
  defined.
  \end{itemize}
In addition
\begin{itemize}
\item $g$ uniformly contracts in the $y$-direction: there
  exists $\mu<1$ such that $|\partial_yg|<\mu$;
\item indeed, $\partial_yg(x,y)\approx x^\beta$ for
  $x\approx0$ where $\beta=-\lambda_2/\lambda_1$ and since
  $\beta>1>\alpha>0$ we have another crucial relation:
  $\lim_{x\to0} \frac{\partial_xg(x,y)}{f'(x)}=0$;
\end{itemize}
this ensures that the contraction in the $y$-direction is
much stronger than the expansion near the singular line
$\ell$, which ensures that the foliation $\cF^s$ is
persistent for all nearby $C^1$ flows. These properties
ensure that $\Lambda$ contains a dense regular orbit and
that this property persists for all $C^1$ close enough
flows, that is, \emph{$\Lambda$ is robustly transitive}.

For a detailed construction of a geometric Lorenz attractor
see~\cite{AraPac2010,galapacif09}.  As mentioned above, a
geometric Lorenz attractor is the simplest example of a
singular hyperbolic attractor~\cite{MPP98}.

Tucker in its PhD thesis \cite{Tu2}
  presented a computer assisted proof that showed the
  existence of a rectangular cross-section $S$ for the flow
  of the Lorenz equations whose first return map $F$ to this
  cross-section $S$ satisfies the same conditions outlined
  above for the geometric Lorenz flow:
\begin{itemize}
\item the image of the $S$ (except the singular
  line) is contained in the interior of $S$,
  which shows there exists an attracting set for the flow
  crossing $S$;
\item on $S$ there exists a field $C^u$ of cones which
  are invariant under the derivative of $F$ and whose
  vectors are expanded by a uniform rate after finitely
  many iterates of $F$ (the algorithm showed that $29$
  iterations are enough to ensure this expansion);
\item the fact that the volume is contracting (the
  divergence of the flow of the Lorenz equations is
  constant and negative) together with the previous
  item, ensures that $F$ contracts area in $S$ and thus
  the field of cones $C^s$, given by the complement of
  the field $C^u$, is invariant for the inverse of $F$
  restricted to the attracting set. This ensures that
  there exists a contracting foliation just like
  $\cF^s$ on $S$ (this is a two-dimensional argument);
\item the one-dimensional quotient map over this
  foliation is piecewise expanding as stated for the
  geometrical Lorenz case.\footnote{Note that we may
    not have that the inverse of the derivative of the
    one dimensional map is of bounded variation. This
    in general will be piecewise Holder and of
    generalized bounded variation as explained in what
    follows. }
\end{itemize}
Hence the attractor of the flow of the Lorenz equations
satisfies the axiomatic conditions stated in the
beginning of this section and, after a suitable
non-linear change of coordinates that rectifies the
contracting foliation $\cF^s$ on $S$, we obtain a flow
similar to a geometric Lorenz attractor (the concrete example described above).

We stress that, in general, two such flows are neither
topologically conjugate nor topologically equivalent;
see for example \cite{GW79} and \cite[Chapter 3,
Section 3]{AraPac2010}.

\section{Singular Hyperbolic systems}
\label{sec:singul-hyperb-system}

In this section we define and describe singular
hyperbolic attractors.  We outline a construction (a suitable Poincar\'e return map) which
allows to rigorously investigate several properties of
its dynamics.  We include a list of properies of the
construction and of the induced map, which are useful
to obtain several results on the statistical properties
of the dynamics, but can be also useful to the reader
for future applications.  
The properties we obtain, shows that the dynamics of this general class of flows is not very different from the ones of the geometric Lorenz one, and thus with some effort is possible to generalize some result obtained for this simple model to the general ones.
The idea is to consider a
suitable Poincar\'e section and a suitable induced
return map.  It turns out that if both are well chosen,
the properties of the return map are similar to the
ones of the geometric Lorenz system described
before. Then by suitable extension of the ideas which
have been used for the geometric case, it is possible
to treat this more general case.

 We say that an attracting set
$\Lambda=\Lambda_X(U)$ for a $3$-flow $X$ and some open
subset $U$ is \emph{robust} if there exists a $C^1$
neighborhood $\cU$ of $X$ in $\fX^1(M)$ such that
$\Lambda_Y(U)$ is transitive for every $Y\in\cU$. 

The following result obtained by Morales, Pacifico and
Pujals in~\cite{MPP04} characterizes robust attractors for
three-dimensional flows.
\begin{theorem}
  \label{thm:robust-attractor-sing-hyp}
  Robust attractors for flows containing equilibria are
  singular-hyperbolic sets.
\end{theorem}
We remark that vector fields exhibiting
robust attractors cannot be $C^1$  approximated by vector fields
presenting either attracting or repelling periodic points
in a neighborhood of the attractor.
This implies that, on $3$-manifolds, any periodic orbit
inside a robust attractor is hyperbolic of saddle-type.

We now precisely define the concept of
singular-hyperbolicity.  A compact invariant set $\Lambda$
of $X$ is {\em partially hyperbolic} if there are a
continuous invariant tangent bundle decomposition $T_\Lambda
M=E^s_\Lambda\oplus E^c_\Lambda$ and constants $\lambda,K>0$
such that
\begin{itemize}
\item
{\em $E^c_\Lambda$ $(K,\lambda)$-dominates $E^s_\Lambda$},
i.e. for all $x\in\Lambda$ and for all $t\ge0$

\begin{align}\label{eq.domination}
  \| DX^t(x)\mid E^s_x\| \leq \frac{e^{-\lambda t}}{K}\cdot
  m(DX^t(x) \mid E^c_x);
\end{align}
where $m(L)$ for a linear map
  $L:(E,|\cdot|_E)\to (F,|\cdot|_F)$ between normed vector
  spaces denotes the \emph{conorm} defined as
  $m(L)=\inf\{|L(v)|_F:|v|_E=1\}$.
\item
$E^s_\Lambda$ is $(K,\lambda)$-contracting: $\|DX^t\mid
E^s_x\|\le K e^{-\lambda t}$ for all $x\in\Lambda$ and for all $t\ge0$.
\end{itemize}
For $x\in \Lambda$ and $t\in\RR$ we let $J_t^c(x)$ be
the absolute value of the determinant of the linear map
$DX^t(x)\mid E^c_x:E^c_x\to E^c_{X^t(x)}$.  We say that the
sub-bundle $E^c_\Lambda$ of the partial hyperbolic set
$\Lambda$ is {\em $(K,\lambda)$-volume expanding} if
$$
J_t^c(x)=\big| \det(DX^t\mid E^c_x) \big|\geq Ke^{\lambda t},
$$
for every $x\in \Lambda$ and $t\geq 0$.

We say that a partially hyperbolic set is {\em
  singular-hyperbolic} if its singularities are hyperbolic
and it has volume expanding central direction.

A \emph{singular-hyperbolic attractor} is a
singular-hyperbolic set which is an attractor as well: an
example is the geometric
Lorenz attractor presented in Section~\ref{sec:geoLorenz};
and also the attractor in the Lorenz
equations~\eqref{e-Lorenz-system} as a consequence of the
work \cite{Tu2} of Tucker (the work o Tucker indeed proves the existence of the attractor in that given ODE).

Any equilibrium $\sigma$ of a 
singular-hyperbolic attractor for a vector field $X$ is such
that $DX(\sigma)$ has only real eigenvalues
$\lambda_2\le\lambda_3\le\lambda_1$ satisfying the same
relations as in the Lorenz flow example:
\begin{align}
  \label{eq:lorenz-like}
  \lambda_2<\lambda_3<0<-\lambda_3<\lambda_1,
\end{align}
which we refer to as \emph{Lorenz-like equilibria};
this is proved in \cite{MPP04}. 

We recall that an compact $X^t$-invariant set $\Lambda$ is
\emph{hyperbolic} if the tangent bundle over $\Lambda$
splits $T_\Lambda M=E^s_\Lambda\oplus E^X_\Lambda \oplus
E^u_\Lambda$ into three $DX^t$-invariant subbundles, where
$E^s_\Lambda$ is uniformly contracted, $E^u_\Lambda$ is
uniformly expanded, and $E^X_\Lambda$ is the direction of
the flow at the points of $\Lambda$.  It is known,
see~\cite{MPP04,AraPac2010}, that a partially hyperbolic set
for a three-dimensional flow, with volume expanding central
direction and without equilibria, is hyperbolic. Hence the
notion of singular-hyperbolicity is an extension of the
notion of hyperbolicity.

\subsection{Main properties}
\label{sec:main-ideas-constr}

Now we show how to work on this kind of systems, as mentionned in the introduction to this section we outline the main ideas and features of a construction, showing that there is a suitable section of the system, having a return map which preserves a contracting foliation, and has several other properties in common with the one of  geometric Lorenz systems.
We outline a construction which is useful to obtain the statistical properties we mention in the following sections. The construction has been modified in the literature to obtain slightly different properties (see \cite{AraPac2010}) but the general strategy is the same.

The main idea is to obtain a family of adapted
cross-sections and Poincar\'e maps between them which, under
a suitable choice of coordinates, can be combined together
to obtain a map $F$ which has properties similar to the ones of the return map in the geometric Lorenz systems. More precisely, we will have the following.

\begin{theorem}
  \label{thm:propert-singhyp-attractor}
  For an open and dense subset of $C^2$ vector fields $X$ having
  a singular hyperbolic attractor $\Lambda$ on a
  $3$-manifold, there exists a finite family $\Xi$ of
  cross-sections and a global ($n$-th return) Poincar\'e map
  $R:\Xi_0\to\Xi$, $R(x)=X_{\tau(x)}(x)$ such that

  \begin{enumerate}
  \item the domain $\Xi_0=\Xi\setminus\Gamma$ is the entire
    cross-sections with a family $\Gamma$ of finitely many
    smooth arcs removed and $\tau:\Xi_0\to[\tau_0,+\infty)$
    is a smooth function bounded away from zero by some
    uniform constant $\tau_0>0$.
  \item
    We can choose coordinates on $\Xi$ so that the map $R$
    can be written as $F:\tilde Q\to Q$, $F(x,y)=(T(x),G(x,y))$,
    where $Q=\II\times\II$, $\II=[0,1]$ and
    $\tilde Q=Q\setminus\Gamma_0$ with $\Gamma_0=\cC\times\II$
    and $\cC=\{c_1,\dots,c_n\}\subset\II$ a finite set of
    points.
  \item The map $T:\II\setminus\cC\to\II$ is $C^{1+\alpha}$
    piecewise monotonic with $n+1$ branches defined on the
    connected components of $\II\setminus\cC$ and has a
    finite set of a.c.i.m., $\mu^i_T$. Also $\inf|T'|>1$ where it is
    defined, $1/|T'|$ has universal bounded $p$-variation
    and then $d\mu^i_T/dm$ has bounded $p$-variation.
  \item
    The map $G:\tilde Q\to\II$ preserves and uniformly contracts
    the vertical foliation
    $\cF=\{\{x\}\times\II\}_{x\in\II}$ of $Q$: there
    exists $0<\lambda<1$ such that $
    \dist(G(x,y_1),G(x,y_2)) \leq \lambda \cdot |y_1-y_2|$
    for each $y_1, y_2 \in \II$. In addition, the map $G$
    satisfies $\varsq(G)<\infty$.
  \item
    The map $F$ admits a finite family of  physical probability measures
    $\mu^{i}_{F}$ which are induced by $\mu^i_T$ in a standard
    way. The Poincar\'e time $\tau$ is integrable both with
    respect to each $\mu^{i}_{F}$ and with respect to the
    two-dimensional Lebesgue area measure of $Q$.
  \item
    Moreover if, for all singularities $\sigma\in\Lambda$,
    we have the eigenvalue relation
    $-\lambda_2(\sigma)>\lambda_1(\sigma)$, then the second
    coordinate map $G$ of $F$ has a bounded partial
    derivative with respect to the first coordinate, i.e.,
    there exists $C>0$ such that $|\partial_x G(x,y)|<C$ for
    all $(x,y)\in(\II\setminus\{c_1,\dots, c_n\})\times\II$.
  \end{enumerate}
\end{theorem}


This result shows that the dynamics of singular hyperbolic attractors has several aspects in common with the one of geometric Lorenz attractors. This allows to extend some method of proof, and obtain most of the few rigorous results on the statistical properties of the singular hyperbolic dynamics.

We now give some idea of the general construction that allows to obtain the result.
Then we will show some technical details about the main steps.

The idea is to take a suitable Poincare section, which is made by a family of rectangles.
Near each fixed point of the flow, we will put six rectangles like in figure \ref{fig:singularbox0}. We know that the dynamics near  the flow is similar to the one of the geometric Lorenz system. 
The two dynamics indeed can be identified by a local linearization of the flow near the fixed point.
The linearization can have different regularity properties, according to the arithmetical properties of the eigenvalues associated to the fixed point. We will need a   $ C^2$  linearization. The precise requirements are outlined in theorem \ref{thm:smooth-linear} below. This linearization allows to obtain precise informations on the behavior of the Poincar\'e maps near  the fixed points, which will have (after linearization) the same form as the geometric Lorenz ones. What said allows to obtain information on the dynamics near the fixed points.
Far from them we have an hyperbolic dynamics with  expanding and contracting directions.
We complete the set of sections we need by taking a suitable number of further sections, intercepting the flow near the attractor. We can take a finite number of them by a  compactness argument.
The hyperbolicity of the flow, impolies that if we choose the Poincar\'e map in a way that  the time which is necessary for an orbit to go from a rectangle to the next one is large enough, then the induced map is hyperbolic.
The regularity of the flow, of the linearization, and of the sections reflects on the regularity of the induced Poincar\'e map, and of the preserved contracting foliation.
The presence and regularity of the preserved one dimensional foliation imply that the map has the form given at item 2.

Now we present some of the main steps in the construction, see (\cite{AraPac2010} and \cite{ArGalPac} ) for further details.

\subsubsection{Linearization near the singularities}
Here we prsent the assumption which are needed in the above theorem.
One step which is important to obtain information on the return map, is the control of its behavior near the fixed point. This can be approached by linearization of the system near the fixed point.

We recall that, in general, hyperbolic singularities are
 linearizable by an H\"older homeomorphism
according to the standard Hartman-Grobman
Theorem~\cite{PM82,robinson1999}.

For our construction, we need a smoother linearization.
A result in this direction is provided by Hartman (See \cite[Theorem 12.1, p. 257]{Hartman02}). 
In the absence of
resonances, orbits of the flow in a small neighborhood
$U_\sigma$ of the given equilibrium $\sigma$ are solutions
of the linear system (\ref{eq:LinearLorenz}), modulo a
smooth change of coordinates with $\lambda_2<\lambda_3<0<-\lambda_3<\lambda_1$.

\begin{theorem}
  \label{thm:smooth-linear}
  Let $n\in\ZZ^+$ be given. Then there exists an integer
  $N=N(n)\ge2$ such that: if $\Gamma$ is a real non-singular
  $d\times d$ matrix with eigenvalues
  $\gamma_1,\dots,\gamma_d$ satisfying
  \begin{align}\label{eq:non-resonance}
    \sum_{i=1}^d m_i \gamma_i \neq \gamma_k
    \quad
    \text{for all}
    \quad
    k=1,\dots, d
    \qand
    2\le\sum_{j=1}^d m_j\le N
  \end{align}
  and if $\dot\xi=\Gamma\xi+\Xi(\xi)$ and
  $\dot\zeta=\Gamma\zeta$, where $\xi,\zeta\in\RR^d$ and
  $\Xi$ is of class $C^N$ for small $\|\xi\|$ with
  $\Xi(0)=0, \partial_\xi\Xi(0)=0$; then there exists a
  $C^n$ diffeomorphism $R$ from a neighborhood of $\xi=0$ to
  a neighborhood of $\zeta=0$ such that $R\xi_t
  R^{-1}=\zeta_t$ for all $t\in\RR$ and initial conditions
  for which the flows $\zeta_t$ and $\xi_t$ are defined in
  the corresponding neighborhood of the origin.
\end{theorem}

Hence it is enough for us to choose the eigenvalues
$(\lambda_1,\lambda_2, \lambda_3)\in\RR^3$ of $\sigma$ satisfying a
\emph{finite set of non-resonance relations}
\eqref{eq:non-resonance} for a certain $N=N(2)$ and for each
singularity $\sigma_k$ in $\Lambda$. For this condition
defines an open and dense set in $\RR^3$ and so all small
$C^1$ perturbations $Y$ of the vector field $X$ will have a
singularity whose eigenvalues
$(\lambda_1(Y),\lambda_2(Y),\lambda_3(Y))$ are still in the
$C^2$ linearizing region.

We note that in (\ref{eq:LinearLorenz}) $x_1$ corresponds to
the strong-stable direction at $\sigma$, $x_2$ to the
expanding direction and $x_3$ to the weak-stable direction.

\begin{figure}[ht]
\psfrag{S1}{$\Sigma^{i,+}$}
\psfrag{S2}{$\Sigma^{i,-}$}
\psfrag{S3}{$\Sigma^{o,+}$}
\psfrag{S4}{$\Sigma^{o,-}$}
\psfrag{s}{$\sigma$}
\psfrag{L1}{$\ell^+$}
\psfrag{L2}{$\ell^-$}
\psfrag{z}{{\footnotesize $z$}}
\psfrag{R}{{\footnotesize $R(z)$}}
\psfrag{X1}{{\footnotesize $x_1$}}
\psfrag{X2}{{\footnotesize $x_2$}}
\psfrag{X3}{{\footnotesize $x_3$}}
\includegraphics[height=5cm]{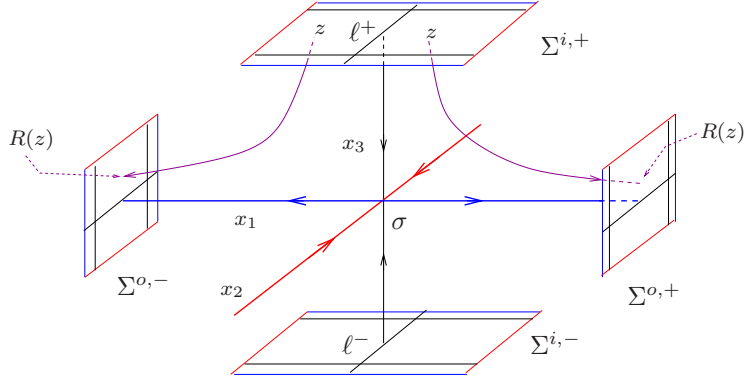}
\caption{\label{fig:singularbox0} Cross-sections near a
  Lorenz-like equilibrium.}
\end{figure}

Then for some $\delta>0$ we may choose
cross-sections contained in $U_\sigma$
\begin{itemize}
\item $\Sigma^{o,\pm}_\sigma$ at points $y^{\pm}$ in
  different components of
  $W^u_{loc}(\sigma)\setminus\{\sigma\}$
\item $\Sigma^{i,\pm}_\sigma$ at points $x^{\pm}$ in
  different components of $W^s_{loc}(\sigma)\setminus
  W^{ss}_{loc}(\sigma)$
\end{itemize}
and Poincar\'e first hitting time maps
$R^\pm:\Sigma^{i,\pm}_\sigma\setminus\ell^\pm\to
\Sigma^{o,-}_\sigma\cup\Sigma^{o,+}_\sigma$, where $
\ell^\pm=\Sigma^{i,\pm}_\sigma\cap W^s_{loc}(\sigma), $
satisfying (see Figure~\ref{fig:singularbox0})
\begin{enumerate}
\item every orbit in the attractor passing through a small
  neighborhood of the equilibrium $\sigma$ intersects some
  of the incoming cross-sections $\Sigma^{i,\pm}_\sigma$;
\item $R^\pm$ maps each connected component of
  $\Sigma^{i,\pm}_\sigma\setminus\ell^\pm$ diffeomorphically
  inside a different outgoing cross-section
  $\Sigma^{o,\pm}_\sigma$, preserving the corresponding
  stable foliations.
\end{enumerate}
Here we write $W^*_{loc}(\sigma), *=s,ss,u$ for the local
invariant stable, strong-stable and unstable manifolds of
the hyperbolic saddle-type singularity $\sigma$ (see
e.g. \cite{PM82}), so that these invariant manifold extend
up to the cross-sections $\Sigma^{i,\pm}$ and
$\Sigma^{o,\pm}$.

We note that at each flow-box near a singularity there are
four cross-sections: two ``ingoing''
$\Sigma_{\sigma}^{i,\pm}$ and two ``outgoing''
$\Sigma_{\sigma}^{o,\pm}$.

Using $C^2$ linearizing coordinates in a flow-box near a
singularity, with the appropriate rescaling, we can assume
without loss of generality that, for a small $\delta>0$,  see
Figure~\ref{fig:singularbox0}
\begin{align*}
  \Sigma^{i,\pm}&=\{ (x_1,x_2,\pm1): |x_1|\le\delta ,
  |x_2|\le\delta\} \quad\text{and}
  \\
  \Sigma^{o,\pm}&=\{ (\pm1,x_2,x_3): |x_2|\le\delta ,
  |x_3|\le\delta\}.
\end{align*}
Then from \eqref{eq:LinearLorenz} we can determine the
expression of the Poincar\'e maps between ingoing and
outgoing cross-sections after the linearization, easily
\begin{equation}
  \label{eq:nonflatsing}
\Sigma^{i,+}\cap\{x_1>0\}\to \Sigma^{0,+},
\quad
(x_1,x_2,1)\mapsto
\big(1,x_2\cdot x_1^{-\lambda_2/\lambda_1},
x_1^{-\lambda_3/\lambda_1}\big).
\end{equation}
The cases corresponding to the other ingoing/outgoing pairs
and signs of $x_1,x_2$ are similar. The possibility to have an explicit form (after change of vaviables) allows to understand several aspects of the regularity of the return map on the section.

\subsubsection{Physical measure and basic properties}

Here we would like to justify the existence of a physical measure for the Poincar\' e map, which in turn implies the existence of the physical measure for the flow.
The starting point is that the induced one dimensional map being piecewise expanding has an absolutely continuous inveriant measure. Outside this, there are only contracting directions, and this measure give rise to a physical measure for the two dimensional Poincar\'e map.

More precisely, since the one dimensional induced map is piecewise expanding, with generalized bounded variation (piecewise Holder) derivative,  from \cite[Lemma 1.4]{Ke85}

\begin{lemma}\label{ftemacim}
  The one-dimensional map $T$ obtained above has finitely
  many ergodic physical measures $\mu_T^1,\dots,\mu_T^l$,
  whose density is a function of $p$-bounded
  variation, and whose ergodic basins
  cover Lebesgue almost all points of $\II$.
\end{lemma}

According to standard constructions described in \cite{APPV}
and \cite[Section 7.3, pp. 225-235]{AraPac2010}, each
physical measure $\mu_T^i$ can be lifted to a physical
measure $\nu_\Lambda^i$ for the flow of $X$ and supported on
the attractor $\Lambda$; more on this in
Proposition~\ref{pr:tau_muF-int} of
Subsection~\ref{sec:integr-global-poinca}. Since a
singular-hyperbolic attractor is transitive, that is, it has
a dense orbit, it follows that there can be only one such
physical measure for the flow in the basin of attraction of
$\Lambda$; see~\cite[Section 7.3.8,
pp. 234-235]{AraPac2010}. Therefore we have (see \cite{APPV} )

\begin{theorem}
  \label{srb}
  Let $\Lambda=\Lambda_X(U)$ be a singular-hyperbolic
  \emph{attractor} of a flow $X\in\fX^2(M)$ on a
  three-dimensional manifold.  Then $\Lambda$ supports a
  unique physical probability measure $\mu$ which is ergodic
  and its ergodic basin covers a full Lebesgue measure
  subset of the topological basin of attraction, i.e.,
  $B(\mu)=W^s(\Lambda)$ Lebesgue mod $0$.  Moreover the
  support of $\mu$ is the whole attractor $\Lambda$.
\end{theorem}

\subsubsection{Integrability of $\tau$, $\log|T^\prime|$ and
  $\log|\partial_yG|$}
\label{sec:integr-global-poinca}

The global Poincar\'e time $\tau$ is
integrable with respect to both the two-dimensional Lebesgue
measure $m$ on $\QQ$ and  the $F$-invariant physical
measure $\mu_F$ on $\QQ$, which lifts to the physical
measure $\nu$ for the flow on the singular-hyperbolic
attractor and itself is a lift of the $T$-invariant
absolutely continuous probability measure $\mu_T$ on $\II$.

\begin{proposition}
  \label{pr:tau_muF-int}
  The global Poincar\'e time $\tau$ is integrable with
  respect to the $F$-invariant physical probability measure
  $\mu_F$ and with respect to $m$.
\end{proposition}

Some other integrability properties will be needed and can
be obtained using the properties of the maps $T$ and $G$;
see \cite{ArGalPac} for more details.

\begin{proposition}
  \label{pr:logTprimeDyG-int}
  We have the following properties:
  \begin{enumerate}
  \item $0<\int \log|T^\prime|\,d\mu_F <\infty$;
  \item $\int -\log|\partial_yG(x,y)|\,d\mu_F <\infty$;
  \item the  maps $y\mapsto \partial_y G(x,y)$ are
    uniformly equicontinuous for
    $x\in\II\setminus\{c_1,\dots,c_n\}$, i.e., outside the
    singularities of the map $T$.
  \end{enumerate}
\end{proposition}

\subsection{Consequences of absence of sinks and sources
  nearby}
\label{sec:absence-sinks-source}

The proof of Theorem~\ref{thm:robust-attractor-sing-hyp}
given in~\cite{MPP04} uses several tools from the theory of
normal hyperbolicity developed first by Ma\~n\'e
in~\cite{Man82} together with the low dimension of the
flow. 

Lorenz-like equilibria are the only ones contained in robust
attractors naturally, since they are the only kind of
equilibria in a $3$-flow which cannot be perturbed into
saddle-connections which generate sinks or sources when
unfolded. We note that since this kind of
  hyperbolic fixed points $\sigma$ for vector fields $X$ has
  only real eigenvalues and a negative real eigenvalue
  $\lambda_2$ strictly smaller than the rest of the spectrum
  of the tangent map $DX(\sigma)$, then it is well-known
  that there exists an invariant strong-stable manifold
  through the fixed point and tangent to the one-dimensional
  eigenspace corresponding to $\lambda_2$.

\begin{proposition}
\label{tlorenzlike}
Let $\Lambda$ be a robustly transitive set of $X \in
{\fX}^1(M)$.  Then, either for $Y=X$ or $Y=-X$, every
singularity $\sigma\in\Lambda$ is Lorenz-like for $Y$ and
satisfies $ W_Y^{ss}(\sigma)\cap\Lambda=\{\sigma\}.
$

\end{proposition}

The following shows in particular that the notion of
singular hyperbolicity is an extension of the notion of
hyperbolicity.

\begin{proposition}
\label{hypinduced}
Let $\Lambda$ be a singular hyperbolic compact set of $X\in
{\fX}^1(M)$.  Then any invariant compact set
$\Gamma\subset \Lambda$ without singularities is uniformly
hyperbolic.
\end{proposition}

A consequence of Proposition \ref{hypinduced} is that every
periodic orbit of a singular hyperbolic set is
hyperbolic. The existence of a periodic orbit in every
singular-hyperbolic attractor was proved recently
in~\cite{BM04} and also a more general result was obtained
in~\cite{AP}.

\begin{proposition}
\label{ccc}
Every singular hyperbolic attractor $\Lambda$ has a dense
subset of periodic orbits.
\end{proposition}

In the same work~\cite{AP} it was announced that every singular
hyperbolic attractor is the homoclinic class associated to
one of its periodic orbits.  Recall that the
\emph{homoclinic class} of a periodic orbit $\cO$ for $X$ is
the closure of the set of transversal intersection points of
it stable and unstable manifold:
$H(\cO)=\overline{W^u(\cO)\pitchfork W^s(\cO)}$. This result
is well known for the elementary dynamical pieces of
uniformly hyperbolic attractors. Moreover, in particular,
the geometric Lorenz attractor is a homoclinic class as
proved in~\cite{B04}. A proof of this property for every
singular-hyperbolic attractor is presented
in~\cite{AraPac2010}.

\begin{proposition}
  \label{pr:homclass}
  A singular-hyperbolic attractor $\Lambda$ for a
  three-dimensional vector field $X$ is the homoclinic class
  $H(\cO)$ of a hyperbolic periodic orbit $\cO$ of
  $\Lambda$.
\end{proposition}


\section{ More on the Ergodic Theory of singular-hyperbolic attractors}
\label{sec:ergodic-theory-singu}

The ergodic theory of singular-hyperbolic attractors is
incomplete.  Many results still are proved only in the
particular case of  (axiomatic) geometric Lorenz flows (sometimes with
additional assumptions) and several automatically extend to
the original Lorenz flow after the work of Tucker
\cite{Tu2}, but demand an extra effort to encompass the full
singular-hyperbolic setting.

We note that a singular-hyperbolic attractor in general
contains finitely many hyperbolic singularities and does not
admit a single connected cross-section which is crossed by
all orbits except the singularities, as is the case of the
geometrical Lorenz attractor and the attractor of the Lorenz
system of equations. 

\subsection{The physical measure is a  $u$-Gibbs state }
\label{sec:robustn-singul-hyper}

It follows from the proof of Theorem \ref{srb}, in~\cite{APPV} that the
singular-hyperbolic \emph{attracting set} $\Lambda_Y(U)$ for
all $Y\in\fX^2(M)$ which are $C^1$-close enough to $X$
\emph{admits finitely many physical measures whose ergodic
  basins cover $U$ except for a zero volume
  subset}. We note that a
  singular-hyperbolic attractor is not necessarily robustly
  transitive: examples of this behavior are known; see
  \cite[Example 5.7]{AraPac2010}.

  Theorem~\ref{srb} shows that typical orbits in the basin
  of every singular-hyperbolic attractor, for a $C^2$ flow
  $X$ on a $3$-manifold, have well-defined statistical
  behavior, i.e. for Lebesgue almost every point the forward
  Birkhoff time average converges, and it is given by a
  certain physical probability measure $\mu$. It was also
  obtained that this measure is hyperbolic and admits
  absolutely continuous conditional measures along the
  center-unstable directions on the attractor. As a
  consequence, it is a $u$-Gibbs state and an equilibrium
  state for the flow.

Here hyperbolicity of the invariant measure $\mu$ means
\emph{non-uniform hyperbolicity} of the probability measure
$\mu$: the tangent bundle over $\Lambda$ splits into a sum
$T_z M = E^s_z\oplus E^X_z\oplus F_z$ of three
one-dimensional invariant subspaces defined for $\mu$-a.e.
$z\in \Lambda$ and depending measurably on the base point
$z$, where $\mu$ is the physical measure in the statement of
Theorem~\ref{srb}, $E^X_z$ is the flow direction (with zero
Lyapunov exponent) and $F_z$ is the direction with positive
Lyapunov exponent, that is, for every non-zero vector $v\in
F_z$ we have
\[
\lim_{t\to+\infty}\frac1t\log\|DX^t(z)\cdot v\|>0.
\]
We note that the invariance of the splitting implies that
$E^{cu}_z=E^X_z\oplus F_z$ whenever $F_z$ is defined.

Theorem~\ref{srb} is another statement of sensitiveness,
this time applying to the whole essentially open set
$B(\Lambda)$.  Indeed, since non-zero Lyapunov exponents
express that the orbits of infinitesimally close-by points
tend to move apart from each other, this theorem means that
most orbits in the basin of attraction separate under
forward iteration.  See Kifer~\cite{Ki88}, and
Metzger~\cite{mtz001}, and references therein, for previous
results about invariant measures and stochastic stability of
the geometric Lorenz models.

The $u$-Gibbs property of $\mu$ is stated as follows. 

\begin{theorem}
\label{thm:srbmesmo}
Let $\Lambda$ be a singular-hyperbolic attractor for a $C^2$
three-dimen\-sional flow.  Then the physical measure $\mu$
supported in $\Lambda$ has a disintegration into absolutely
continuous conditional measures $\mu_\gamma$ along
center-unstable surfaces $\gamma$ such that
$\frac{d\mu_\gamma}{dm_\gamma}$ is uniformly bounded from
above.  Moreover $\supp(\mu)=\Lambda\,$.
\end{theorem}

\subsection{Entropy formula}

Here the existence of unstable manifolds is guaranteed by
the hyperbolicity of the physical measure: the
strong-unstable manifolds $W^{uu}(z)$ are the ``integral
manifolds'' in the direction of the one-dimensional
sub-bundle $F$, tangent to $F_z$ at almost every
$z\in\Lambda$. The sets $W^{uu}(z)$ are embedded
sub-manifolds in a neighborhood of $z$ which, in general,
depend only measurably (including its size) on the base
point $z\in\Lambda$.

We remark that since $\Lambda$ is an
attracting set, then $W^{uu}(z)\subset\Lambda$ whenever
defined. The central unstable surfaces mentioned in the
statement of Theorem~\ref{thm:srbmesmo} are just small 
strong-unstable manifolds carried by the flow, which are
tangent to the central-unstable direction $E^{cu}$.

The absolute continuity property along the center-unstable
sub-bundle given by Theorem~\ref{thm:srbmesmo} ensures that
\[
h_\mu(X^1)=\int \log\big| \det (DX^1\mid E^{cu})  \big| \, d\mu,
\]
by the characterization of probability measures satisfying
the Entropy Formula, obtained in~\cite{LY85}.  The above
integral is the sum of the positive Lyapunov exponents along
the sub-bundle $E^{cu}$ by Oseledets
Theorem~\cite{Man87,Wa82}.  Since in the direction $E^{cu}$
there is only one positive Lyapunov exponent along the
one-dimensional direction $F_z$, $\mu$-a.e. $z$, the
ergodicity of $\mu$ then shows that the following is true.

\begin{corollary}
    \label{cor:uGibbs}
  If $\Lambda$ is a singular-hyperbolic attractor for a
  $C^2$ three-dimen\-sional flow $X^t$, then the physical
  measure $\mu$ supported in $\Lambda$ satisfies the Entropy
  Formula
\[
h_\mu(X^1)=\int\log\| DX^1\mid F_z\|\, d\mu(z).
\]
\end{corollary}
From the characterization of measures satisfying the Entropy
Formula given in \cite{LY85}, we see that \emph{$\mu$ has
  absolutely continuous disintegration along the
  strong-unstable direction}, along which the Lyapunov
exponent is positive, thus \emph{$\mu$ is a $u$-Gibbs
  state}~\cite{PS82}.  This also shows that \emph{$\mu$ is
  an equilibrium state for the potential} $-\log\| DX^1\mid
F_z\|$ with respect to the diffeomorphism $X^1$. We note
that the entropy $h_\mu(X^1)$ of $X^1$ is the entropy of the
flow $X^t$ with respect to the measure $\mu$ \cite{Wa82}.

Hence we are able to extend a basic result on the ergodic
theory of smooth hyperbolic attractors to the setting of
smooth singular-hyperbolic attractors: a hyperbolic
attractor of a $C^2$ diffeomorphism admits a physical
measure which is the only equilibrium state with respect to
the potencial $\log|\det(Df\mid E^u)|$ given by the norm of
the Jacobian of the map along the unstable directions at the
attractor.

\subsection{Exact dimensionality}
\label{sec:muexata}

We recall a result of Steinberger \cite{Stein00} about the local dimension of
Lorenz like systems and prove that for the singular
hyperbolic system the local dimension is defined at almost
every point.

Let us consider a map $F:\QQ\to \QQ $,
$F(x,y)=(T(x),G(x,y))$ where
\begin{enumerate}
\item[(1)] $T:[0,1]\to [0,1]$ is piecewise monotonic: there
  are $c_i\in [0,1]$ for $0\leq i\leq N$ with $0=c_0< \cdots
  < c_N=1$ such that $T|(c_i,c_{I+1})$ is continuous and
  monotone for $0\leq i < N$. Furthermore, for $0\leq i< N$,
  $T|(c_i,c_{i+1})$ is $C^1$ and that $\inf_{x\in
    \cP}|T'(x)|> 0$
  holds where $\cP=[0,1]\setminus \cup_{0\leq i < N}c_i$.
\item[(2)] $G:\QQ\to (0,1)$ is $C^1$ on $\cP\times [0,1]$. Furthermore,
$\sup|\partial G/\partial x| < \infty$, $\sup|\partial G/\partial y| < 1$ and
$|(\partial G/\partial y)(x,y)| > 0$ for $(x,y)\in \cP\times [0,1]$.
\item[(3)] $F((c_i,c_{i+1})\times [0,1])\cap
  F((c_j,c_{j+1})\times [0,1])=\emptyset$ for distinct $i,
  j$ with $0\leq i, j < N$.
\end{enumerate}
Now consider the projection $\pi_{x}:\QQ\to I$, set
${\cV}=\{(c_i,c_{i+1}), 1\leq i \leq N \}$ and
$\cV_k=\bigvee_{i=0}^{k}f^{-i}\cV$.  For $x\in E$ let
$J_k(x)$ be the unique element of $\cV_k$ which contains
$x$.  We say that $\cV$ {\em is a generator} if the length
of the intervals $J_k(x)$ tends to zero for $n \to \infty$
for any given $x$. In piecewise expanding maps it is easy to
see that ${\cV}$ is a generator.  Set
\begin{align}\label{eq:log-psi-vfi}
  \psi(x,y)=\log|T'(x)| \quad \mbox{and} \quad
  \varphi(x,y)=-\log|(\partial G/\partial y)(x,y)|.
\end{align}
The result of Steinberger that we shall use is the following

\begin{theorem}\cite[Theorem 1]{Stein00}
  \label{th:Steinberg} Let $F$ be a two-dimensional map as
  above and $\mu$ an ergodic, $F$-invariant probability
  measure on $I^2$ with the entropy $h_{\mu}(F)>0$.  Suppose
  $\cV$ is a generator, $\int\varphi\cdot d\mu_F < \infty$
  and $0<\int \psi d\mu_F < \infty$.  If the maps $y\mapsto
  \varphi(x,y)$ are uniformly equicontinuous for $x \in
  I\setminus\{0\}$ and $1/|T'|$ has finite universal $p$- Bounded
  Variation, then
$$
d_{\mu}(x,y)=h_{\mu}(F)\big(\frac{1}{\int\psi\cdot d\mu} +
\frac{1}{\int \varphi\cdot d\mu}\big)
$$
for $\mu$-almost all $(x,y) \in \QQ$.
\end{theorem}

Item (3) above is satisfied in our case because the map is
induced by a first return Poincar\'e map induced by a flow.  Moreover
$\sup|\partial G/\partial x| < \infty$ in item (2) above is
established at item (6) of Theorem
\ref{thm:propert-singhyp-attractor}, provided that for all
equilibria $\sigma\in\Lambda$ we have the eigenvalue
relation $-\lambda_2(\sigma)>\lambda_1(\sigma)$.  

Let us also observe that, for the first return map
$F:\QQ\setminus \Gamma\to \QQ$, associated to the
singular-hyperbolic flow, the entropy is positive $h_\mu(F)
> 0$. Indeed, since we know that $\pi\circ F=T\circ\pi$, where
$\pi:\QQ\to\II$ is the projection on the first coordinate,
and that $h_{\mu_T}(T)=\int\log|DT|\,d\mu_T>0$, where
$\mu_T$ is the unique ergodic absolutely continuous
$T$-invariant probability measure, we see that
$h_{\mu_F}(F)>0$.

So, all we need to prove that $(\Xi,F,d\mu_F)$ is exact
dimensional is to verify that $F(x,y)$ satisfies the
hypothesis of Theorem \ref{th:Steinberg}, where
$F:\Xi\circlearrowleft$ is the Poincar\'e return map to
the family of cross-sections $\Xi$ described in
Section~\ref{sec:singul-hyperb-system}; and $\mu_F$ is
the $F$-invariant ergodic SRB measure induced on $\Xi$
by the physical measure of the attractor.

However, Proposition~\ref{pr:logTprimeDyG-int} provides
precisely that for the functions $\vfi,\psi$ defined
above in (\ref{eq:log-psi-vfi}): we have
\begin{enumerate}
 \item $\int \varphi d\mu_F < \infty$;
\item $0< \int \psi d\mu_F < \infty$; and
\item the maps $y \mapsto \varphi(x,y)$ are uniformly equicontinuous
for $x\in\II\setminus\{c_1,\dots,c_n\}$.
\end{enumerate}

This all together finishes the proof of
Theorem~\ref{th:Steinberg} establishing that $\mu_F$ is
exact dimensional.

The exact dimensionality of the measure on the section
implies the exact dimensionality of the measure $\mu$ on the
flow at almost each point, and the dimension satisfies
$d_{\mu}(x)=d_{\mu_F}(x)+1$ at almost every point $x$.


\subsection{Decay of correlations for Poincar\'e maps on
  singular-hyperbolic attractors}
\label{sec:decay-correl-poincar-1}

After obtaining an interesting invariant probability measure
for a dynamical system the next thing to do is to study the
properties of this measure. Besides ergodicity there are
various degrees of mixing (see
e.g. \cite{Wa82,Man87}).

\subsubsection{Decay of correlations for maps versus flows}
\label{sec:Hypflows-decay}


Given a flow $X$ and an invariant ergodic probability
measure $\mu$, we say that the system $(X,\mu)$ is
\emph{mixing} if for any two measurable sets $A,B$
\begin{align}
  \label{eq:defmixA}
\mu\big(A\cap X^{-t}B\big)\xrightarrow[t\to\infty]{}
\mu(A)\cdot\mu(B)
\end{align}
or equivalently
\begin{align*}
\int\vfi\cdot\big(\psi\circ X^t\big) \, d\mu
\xrightarrow[t\to\infty]{}
\int\vfi\,d\mu\int \psi\, d\mu
\end{align*}
for any pair $\vfi,\psi:M\to\RR$ of continuous functions.

Considering $\vfi$ and $\psi\circ X^t:M\to\RR$ as random variables
over the probability space $(M,\mu)$, this definition just
says that ``the random variables $\vfi$ and $\psi\circ X^t$ are
asymptotically independent'' since the expected value
$\EE\big(\vfi\cdot(\psi\circ X^t)\big)$ tends to the product
$\EE(\vfi)\cdot\EE(\psi)$ when $t$ goes to infinity. The
\emph{correlation function}
\begin{align}
  \label{eq:defcorrfunction}
C_t(\vfi,\psi)
&=
\big|\EE\big(\vfi\cdot(\psi\circ X^t)\big)
- \EE(\vfi)\cdot\EE(\psi)\big|\nonumber
\\
&=
\Big|\int\vfi\cdot\big(\psi\circ X^t\big) \, d\mu
-\int\vfi\,d\mu\int \psi\, d\mu\Big|
\end{align}
satisfies  $C_t(\vfi,\psi) 
\xrightarrow[t\to\infty]{}0$ in this case. The \emph{rate of approach to
  zero of the correlation function} is called \emph{the rate
  of decay of correlations} for the observables $\vfi$ and $\psi$
of the system $(X,\mu)$.

The study of decay of correlations for hyperbolic systems
goes back to the work of Sinai~\cite{Si72} and
Ruelle~\cite{Ru76}. Many results were obtained for
transformations. For a diffeomorphism $f$ the notion of decay
of correlations is the same as above replacing $X^t$ by
$f^n$ and letting $n$ go to infinity. Since \cite{Bo75,Ru76}
it is known that the \emph{physical} (SRB) measures for
Axiom A \emph{diffeomorphisms} are mixing and have
\emph{exponential decay of correlations}, that is, there
exists a constant $\alpha\in(0,1)$ such that given
$\vfi$ and $\psi$ there exists $C=C(\vfi,\psi)>0$ such that
\begin{align}
  \label{eq:defexpdecaycorrelations}
  C_n(\vfi,\psi) \le C\cdot e^{-\alpha n} \quad\text{for
    all}\quad n\ge1,
\end{align}
for a suitable class of continuous functions $M\to\RR$, in
this case the H\"older continuous functions.

In more general cases for smooth endomorphisms (see
e.g. \cite{Ho03,alves-luzzatto-pinheiro2005} and references
therein) where the inverse in \eqref{eq:defmixA} is to be
taken as the inverse image of $f^n$, it is possible to have
slower rates of decay.

In contrast to the results available in the case of discrete
dynamical systems, obtaining the rate of decay of
correlations for flows seems to be much more complex and
some results have been established for Anosov flows only
recently.  Ergodicity and mixing for geodesic flows on
manifolds of negative curvature are known since the early
half of the XXth century \cite{hopf39,AS67,Sinai60}.

The proof of exponential decay of correlations for geodesic
flows on manifolds of constant negative curvature was first
obtained in two \cite{CEG84,moore87,ratner87} and three
dimensions \cite{pollicott92} through group theoretical
arguments.

\subsubsection{Decay of correlations for fiber contracting
  maps}
\label{sec:decay-correl-fiber}

In  \cite{ArGalPac} we establish results
on the decay of correlations and convergence to equilibrium
for fiber contracting maps with a fastly converging to
equilibrium base, which imply the following statement (and
this in turn is then applied to singular hyperbolic
attractors).

We recall
  that a measurable map $h:[a,b]\to\RR$ is of \emph{universal
  $p$-bounded variation} if
  \begin{align*}
    \sup_{a=a_0<a_1<\dots<a_n=b} \left( \sum_{j=1}^n \big|
      h(a_i)-h(a_{i-1}) \big|^{1/\alpha} \right)^\alpha
    <\infty,
  \end{align*}
  where the supremum is taken over all finite partitions of
  the interval $I=[a,b]$.

  We will need another definition of variation for maps with
  two variables.  Similarly to the one dimensional case, if
  $f:Q\rightarrow \mathbb{R}$ and $x_{i}\leq x_{2}\leq
  ...\leq x_{n}$, let us define
\begin{equation*}
var^{\square }(f,x_{1},...,x_{n},y_{1},...,y_{n})=\sum_{1\le
  i\leq
n}|f(x_{i},y_{i})-f(x_{i+1},y_{i})|.
\end{equation*}
We then consider the supremum $var^{\square
}(f,x_{1},...,x_{n},y_{1},...,y_{n})$ over all subdivisions $x_{i}$ and all
choices of the $y_{i}$
\begin{equation*}
  var^{\square }(f)=\sup_{n}\left( \sup_{(x_{i}\leq x_{2}\leq ...\leq
      x_{n})\in\II,(y_{i})\in\II}var^{\square
    }(f,x_{1},...,x_{n},y_{1},...,y_{n})\right).
\end{equation*}

  We recall that we denote by $Q=\II\times\II$ the unit
  square, where $\II=[0,1]$. For a function $g:Q\to\RR$ we
  denote by $L(g)$ the best Lipschitz constant of $g$, that
  is, $L(g)=\sup_{p,q\in Q}\frac{|g(p)-g(q)|}{|p-q|}$ where
  $|\cdot|$ is the Euclidean distance. We define the
  Lipschitz norm by setting $\Vert g\Vert_{lip}=\Vert g\Vert
  _{\infty }+L(g)$ where, as usual,
  $\|g\|_\infty=\esssup_{p\in Q}|g(p)|$ and set
  $\textrm{Lip}(Q)=\{g:Q\to\RR:\|g\|_{lip}<\infty\}$.

\begin{theorem}\label{mthm:expdecay}
  Let us consider a map $F:Q\circlearrowleft$ from the unit
  square into itself such that:
\begin{enumerate}
\item $F$ has the form $F(x,y)=(T(x),G(x,y))$ (is a
  skew-product and preserves the natural vertical foliation
  of the square) ;
\item $F|_{\gamma }$ is $\lambda $-Lipschitz with $\lambda
  <1$ (hence is uniformly contracting) on each leaf $\gamma
  $ of the vertical foliation of the square;
\item $var^{\square }(G)<\infty$;
\item $T:\II\circlearrowleft$ is piecewise monotonic, with
  $n+1$, $C^{1}$ increasing branches on the intervals
  $(0,c_{1})$,...,$%
  (c_{i},c_{i+1})$ ,..., $(c_{n},1)$ and $\inf_{x\in\II}
  |T^{\prime }(x)|>1$.
\item $\frac{1}{T^{\prime }}$ has finite universal
  $p-$bounded variation (as defined above);
\item $T$ has only one absolutely continuous
  (w.r.t. Lebesgue on $\II$) invariant probability measure
  (a.c.i.m.) for which it is weakly mixing.
\end{enumerate}
Then {\em the unique physical measure} $\mu _{F}$ of $F$ has
exponential decay of correlation with respect to Lipschitz
observables, that is, there are $C,\Lambda \in
\mathbb{R}^{+},~\Lambda <1$, such that
\begin{align*}
  \left|\int f\cdot(g\circ F^{n})\,d\mu _{F} -
    \int g~d\mu _{F} \int f~d\mu _{F} \right|\leq
  C\Lambda^{n}
  \|g\|_{Lip}\|f\|_{Lip}, \quad f,g\in\textrm{Lip}(Q).
\end{align*}
\end{theorem}

We remark that items (4) to (6) of the assumptions on the
above theorem can be replaced by (more general) exponential
convergence to equilibrium on the base map under suitable
observables. 

\subsubsection{Decay of correlations for the Poincar\'e
  return map of singular-hyperbolic attractors}
\label{sec:decay-correl-poincar}

We now apply these results to singular-hyperbolic attractors
for three-dimensional flows. We let $SH^2(M^3)$ be the
family of all $C^2$ vector fields $X$ on a compact
three-manifold having an open trapping region $U$, i.e.,
$\overline{X^t(U)}\subset U$ for all $t>0$, such that its
maximal invariant subset $\Lambda=\cap_{t>0}X^t(U)$ is a
compact transitive singular-hyperbolic set. We consider the
$C^2$ topology of vector fields in $SH^2(M^3)$ in what follows.

The proof of the existence of physical measure in
Theorem~\ref{thm:srbmesmo} and Corollary~\ref{cor:uGibbs} is
based on a construction that we describe in the next
subsection, and which suitably explored yields the
assumptions on Theorem~\ref{mthm:expdecay}.

This construction shows that there exists a finite family
$\Xi$ of well-adapted cross-sections of the flow on the
attractor where we can define a Poincar\'e return map $F$
which satisfies the properties in the statement of
Theorem~\ref{mthm:expdecay} after a suitable choice of
coordinates.  We remark that, to take advantage of a result
from \cite{Stein00} on exact dimensionality of certain
classes of measures, we need that the Poincar\'e return map
$F$ be injective, which is not evident in the construction
and choice of these cross-sections at \cite{APPV}. Because
of this, the construction presented here is slightly
different; see Subsection \ref{sec:main-ideas-constr} for
details.
 
We need some conditions on the eigenvalues of
the equilibria of $X$ inside $\Lambda$ to reduce the
dynamics on $\Lambda$ to a map $F$ as in
Theorem~\ref{mthm:expdecay}, to obtain
\begin{corollary}
  \label{cor:decaycorr-sing-hyp}
  There exists an open dense set $\cA$ of vector fields in
  $SH^2(M^3)$ such that, for each $X\in\cA$, we can find a
  finite family $\Xi$ of cross-sections to the flow $X_t$ of
  $X$ whose Poincar\'e first return map
  $F:\textrm{dom}(F)\subset\Xi\to\Xi$ has a unique SRB measure $\mu
  _{F}$ which has exponential decay of correlations with
  respect to Lipschitz observables: there are $C,\Lambda \in
  \mathbb{R}^{+},~\Lambda <1$ satisfying for every pair
  $f,g:\Xi\to\RR$ of Lipschitz functions
\begin{equation*}
  \left|\int f\cdot(g\circ F^{n})~d\mu -\int g~d\mu \int
    f~d\mu \right|
  \leq C\Lambda
  ^{n}||g||_{Lip}||f||_{Lip}, \quad n\ge1.
\end{equation*}
\end{corollary}

We remark that, since the works of
Ruelle~\cite{ruelle1983} and
Pollicott~\cite{pollicott99} it is well-known that
exponentially mixing for a base transformation of a
suspension flow does not imply fast mixing for the
suspension flow. In fact, in the suspension flow can be
non-mixing!  Hence we cannot deduce any kind of mixing
results for the flow on a singular-hyperbolic attractor
from Corollary~\ref{cor:decaycorr-sing-hyp}. However,
in a recent work of one of the authors with Varandas
\cite{ArVar}, described at
Section~\ref{sec:recent-results}, it has been proved
the existence of a $C^2$ open subset of vector fields
having a geometrical Lorenz attractor with exponential
decay of correlations for the flow on $C^1$
observables, from which follows the exponential decay
of correlations for the corresponding Poincar\'e
map. But this $C^2$ open subset was obtained under very
strong conditions which cannot hold in such generality
as in Corollary~\ref{cor:decaycorr-sing-hyp}. In
another recent work \cite{AMV13} by the same authors
together with Melbourne, it was proved that \emph{all
  $C^\infty$ geometric Lorenz attractors} (including
the attractor for the system of
equations~\eqref{e-Lorenz-system}) have superpolynomial
decay of correlations, which provides rapid decay of
correlations for the corresponding Poincar\'e map, but
still slower than exponential.

\subsection{Logarithm law for singular hyperbolic attractors}
\label{sec:logarithm-law-syngul}

The decay of correlation for the return map on the
section and the exact dimensionality, implie an
estimation for the behavior of hitting times of
shrinking targets which is called Logarithm law.  This
result essentially says that the time needed to hit a
small target scales as the inverse of the measure of
the target.

Let us recall the a result about discrete time systems we will use: let $(X,F,\mu )$ be an ergodic, measure
preserving transformation on a metric space $X$. Let us
consider a family of target sets $S_{r}$ indexed by a real
parameter $r$ and the time needed for the orbit of a point
$x$ to enter in $S_{r}$
\begin{equation*}
\tau _{F}(x,S_{r})=\min \{n\in \mathbb{N}^{+}:F^{n}(x)\in S_{r}\}.
\end{equation*}
We consider target sets of the form: $S_{r}=\{x\in
X~,~f(x)\leq r\}$, where $f:X\rightarrow \mathbb{R}^{+}$ is
a Lipschitz function, together with the limits
\begin{equation}
\overline{d}(f)=\underset{r\rightarrow 0}{\lim \sup }\frac{\log \mu (S_{r})}{%
\log (r)}~,~\underline{d}(f)=\underset{r\rightarrow 0}{\lim \inf }\frac{\log
\mu (S_{r})}{\log (r)}
\end{equation}
representing a sort of local dimension (the formula for the
local dimension of $\mu $ at a point $x_{0}$ is obtained
when $f(x)=d(x,x_{0})$). When the above limits coincide, we
set $d(f)=\underline{d}(f)=\overline{d}(f).$ In this
setting, the following result is proved in
\cite{galatolo10}; see also \cite{galat07}.
\begin{proposition}
  \label{maine}Let $f$ and $S_{r}$ be as above. Then for
  $\mu$-almost every $x$
\begin{equation}\label{easy}
  \lim \sup_{r\rightarrow 0}\frac{\log \tau _{F}(x,S_{r})}{-\log r}\geq
  \overline{d}(f)~,~\lim \inf_{r\rightarrow 0}\frac{\log \tau _{F}(x,S_{r})}{
    -\log r}\geq \underline{d}(f).
\end{equation}
Moreover, if the system has super-polynomial decay of
correlations under Lipschitz observables and $d(f)$ exists,
then for $\mu$-almost every $x$ it holds
\begin{equation}\label{log-law-disc-time}
\lim_{r\rightarrow 0}\frac{\log \tau _{F}(x,S_{r})}{-\log r}=d(f).
\end{equation}
\end{proposition}

\begin{remark}\label{log-sec}
  Since we are dealing with a ratio of logarithms, and
  \ref{easy} always hold, if we establish the Logarithm Law
  \ref{log-law-disc-time} for some iterate $F^n$, then it
  will hold also for $F$. 
\end{remark}

Let us now see how to extend  the result to flows. Let $X$ be a metric
space, $\Phi ^{t}$ be a measure preserving flow and $\Sigma$
be a section of $(X,\Phi ^{t})$. If the
flow is ergodic and the return time is integrable, then the
hitting time scaling behavior of the flow can be estimated
by the one of the system induced on the section. Hence we
can have a logarithm law for the flow if we can prove it on
the section (with the induced return map).

Given any $x\in X$ let us denote by $t(x)$ the smallest
strictly positive time such that $\Phi ^{t(x)}(x)\in \Sigma
$. We also consider $t^{\prime }(x)$, the smallest non
negative time such that $\Phi^{t^{\prime }(x)}(x)\in
\Sigma$.  We define $\pi:X\rightarrow \Sigma$ as
$\pi(x)=\Phi ^{t^{\prime }(x)}(x)$, the projection on
$\Sigma $. We also denote by $\mu _{F}$ the invariant
measure for the Poincar\'{e} map $F$ which is induced by the
invariant measure $\mu$ of the flow.

\begin{proposition}[\cite{GalNis11}]
  \label{4log}Let us suppose that the flow $\Phi ^{t}$ is
  ergodic and has a section $\Sigma $ with an induced map
  $F$ and invariant measure $\mu _{F}$ such that
$
\int_{\Sigma }t(x)~d\mu _{F}<\infty .
$
Let $r\geq 0$ and $S_{r}\subseteq \Sigma $ be a decreasing
family of measurable subsets with $\lim_{r\rightarrow 0}\mu
_{F}(S_{r})=0$. Let us consider the hitting time relative to
the Poincar\'{e} map
\begin{equation}
\tau ^{\Sigma }(x,S_{r})=\min \{n\in \mathbb{N}^{+};F^{n}(x)\in S_{r}\}.
\end{equation}

Then, there is a full measure set $C\subseteq X$ such that if $x\in C$%
\begin{eqnarray}
\liminf_{r\rightarrow 0}\frac{\log \tau (x,S_{r})}{-\log r}
&=&\liminf_{r\rightarrow 0}\frac{\log \tau ^{\Sigma }(\pi (x),S_{r})}{-\log r%
}, \\
\limsup_{r\rightarrow 0}\frac{\log \tau (x,S_{r})}{-\log r}
&=&\limsup_{r\rightarrow 0}\frac{\log \tau ^{\Sigma }(\pi (x),S_{r})}{-\log r%
}.
\end{eqnarray}
\end{proposition}

We proved in the previous
subsections that the physical invariant measure is exact
dimensional and has exponential decay of correlations.
 Hence by Proposition \ref{maine}  a logarithm law must hold for an iterate of the first return map on the section.
As remarked above, once we have the logarithm law for the
iterate of the first return map, we obtain it for the first
return map. By Remark
\ref{log-sec} and Proposition \ref{4log} we have the
possibility to extend the logarithm law to the flow.

\begin{theorem}\label{thm:exactdim-loglaw}
  Let $\Phi^t:X\circlearrowleft$ be a flow having a singular
  hyperbolic attractor, and let us consider its physical
  invariant measure $\mu$.  Let us consider $x_0$ and the
  local dimension at $x_0$ (which was above proved to exist)
\begin{equation}
  d_{\mu}(x_0)=\lim_{r\rightarrow 0}\frac{\log \mu (B_{r}(x_0))}{\log r},
\end{equation}
then for $\mu$ almost every $x$
\begin{equation}
\lim_{r\rightarrow 0}\frac{\log \tau (x,B_{r}(x_0))}{-\log r}=d_{\mu}(x_0) -1
\end{equation}
where $ \tau (x,B_{r}(x_0))$ is the time needed for the
orbit of $x$ to hit the ball $B_{r}(x_0)$ as above.
\end{theorem}

Of course, noting that Proposition \ref{maine} holds for
targets which are sublevels of Lipschitz functions, it is
possible to give other statements, with the same methods,
replacing balls with more general shrinking targets.


\section{Ergodic theoretic results for the geometric Lorenz
  attractor}
\label{sec:results-geometr-lore}

\subsection{Large Deviations for Lebesgue measure on a
  neighborhood of a geometric Lorenz flow}
\label{sec:large-deviations}

Having shown that physical probability measures exist, it is
natural to consider the rate of convergence of the time
averages to the space average, measured by the volume of the
subset of points whose time averages stay away from the
space average by a prescribed amount up to some evolution
time.
More precisely, if we set $\epsilon>0$ as
an error margin and consider
\[
B_t=\Big\{
z: \Big|
\frac1t\int_{0}^{t}
  \psi\big( X^t(z) \big) 
-
\int\psi\,d\mu
\Big|>\epsilon
\Big\}
\]
then we search conditions under which the
Lebesgue measure of $B_t$ decays to zero exponentially fast,
i.e. weather there are constants $C,\xi>0$ such that
\begin{align*}
\lambda\big( B_t\big) \le C e^{-\xi t}
\quad\mbox{for all}\quad t>0.
\end{align*}
The values of $C,\xi>0$ above depend on $\epsilon,\psi$ and
on global invariants for dynamics.

An extension of part of the results on large deviation
rates of Kifer~\cite{kifer90} from the hyperbolic
setting to semiflows over non-uniformly expanding base
dynamics and unbounded roof function was obtained in
\cite{araujo2006a}. These special flows model
non-uniformly hyperbolic flows like the flow on
singular-hyperbolic attractors. Related results (in
fact, sharper) where obtained by Melbourne and Nicol
\cite{melnicol07} for suspension flows over Markov
towers \emph{assuming that the roof function is
  bounded} with respect to the physical probability
measure of these systems.

\subsubsection{Suspension semiflows}
\label{sec:suspens-semifl}

We first present these flows, present the general strategy
of the approach and then state the main assumptions related
to the modelling of the geometric Lorenz attractor.

Given a H\"older-$C^1$ local diffeomorphism
$f:M\setminus\cS\to M$ outside a volume zero singular set
$\cS$, we say that $\cS$ is \textbf{non-flat} if $f$ behaves
like a power of the distance to \( \cS \): $\| Df(x) \|
\approx \dist(x,\cS)^{-\be}$ for some $\beta>0$; see
Alves-Araujo~\cite{alves-araujo2004} for a precise
statement.

Let also $X^t:M_r\to M_r$ be a \emph{semiflow with roof
  function $r:M\setminus\cS\to\RR$ over the base
  transformation $f$}, as follows.

\begin{figure}[ht]
  \psfrag{8}{$\infty$}
  \psfrag{W}{$X^t_r(x_0,y_0)$}
  \psfrag{r1}{$-r(x_1)$}
  \psfrag{r2}{$-r(x_2)$}
  \includegraphics[width=10cm,height=4cm]{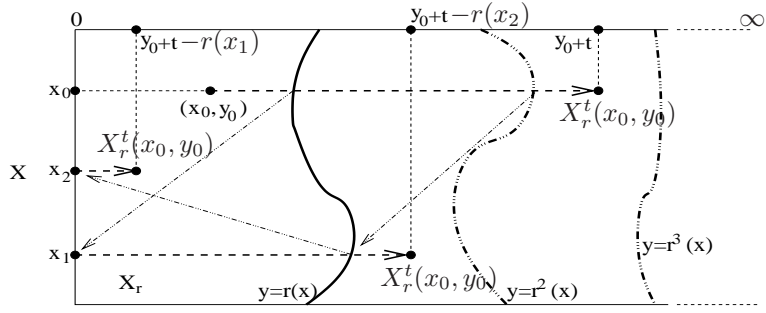}
  \caption{The equivalence relation defining the suspension
    flow of $f$ over the roof function $r$.}
  \label{fig:suspension}
\end{figure}

Where  $M_r=\{ (x,y)\in M\times[0,+\infty): 0\le y < r(x) \}$
and $X^0$ is the identity on $M_r$, where $M$ is a compact
Riemannian manifold.  For $x=x_0\in M$ denote by $x_n$ the
$n$th iterate $f^n(x_0)$ for $n\ge0$.  Denote $S_n^f \vfi(x_0)
= \sum_{j=0}^{n-1} \vfi( x_j )$ for $n\ge1$ and for any
given real function $\vfi$.  Then for each pair
$(x_0,s_0)\in X_r$ and $t>0$ there exists a unique $n\ge1$
such that $S_n^f r(x_0) \le s_0+ t < S_{n+1}^f r(x_0)$ and
define (see Figure~\ref{fig:suspension})
\begin{align*}
   X^t(x_0,s_0) = \big(x_n,s_0+t-S^f_n r(x_0)\big).
\end{align*}

The study of suspension (or special) flows is motivated by
modeling a flow admitting a cross-section. Such flow is
equivalent to a suspension semiflow over the Poincar\'e
return map to the cross-section with roof function given by
the return time function on the cross-section. This is a
main tool in the ergodic theory of hyperbolic flows
developed by Bowen and Ruelle~\cite{BR75}.

The general strategy of the approach is to consider a
suspension semiflow over a piecewise expanding base transformation
with a singular point just like the singularity of the
one-dimensional Lorenz transformation. The result on large
deviations for the Lebesgue measure (note that this measure
is \emph{not} an invariant measure for the flow!) is proved
for this class of suspension flows under suitable conditions
on the base map.

The uniformly contracting foliation on the global
cross-section of a geometric Lorenz flow can be used to show
that the estimates of time averages for points in the
attractor essentially depend only on the stable leaves where
the points lie. This enables a reduction of the estimates of
time averages to time averages for the suspension semiflow
over the one-dimensional Lorenz map. Then it must be checked
that Lorenz one-dimensional map satisfies all the dynamical
assumptions needed to obtain the large deviations bound for
the suspension semiflow.

\subsubsection{Conditions on the base dynamics}
\label{sec:condit-base-dynamics}

We assume that the singular set $\cS$ (containing the points
where $f$ is either \emph{not defined}, \emph{discontinuous}
or \emph{not differentiable}) is regular, e.g. a submanifold
of $M$, and that $f$ is \emph{non-uniformly expanding}:
there exists $c>0$ such that for Lebesgue almost every $x\in
M$
\begin{align*}
  \limsup_{n\to+\infty}\frac1n S_n\psi(x) \le -c
\quad\mbox{where}\quad
\psi(x)=\log\big\| Df(x)^{-1} \big\|.
\end{align*}

Moreover we assume that $f$ has \emph{exponentially slow
  recurrence to the singular set $\cS$} i.e. for all
$\epsilon>0$ there is $\delta>0$ s.t.
\begin{align}\label{eq:expslowrec}
  \limsup_{n\to+\infty}\frac1n\log
\leb\left\{x\in M:
    \frac1nS_n\big| \log d_{\delta}(x,\cS) \big|
    >\epsilon\right\}<0,
\end{align}
where $d_\delta(x,y)=\dist(x,y)$ if $\dist(x,y)<\delta$ and
$d_\delta(x,y)=1$ otherwise.

These conditions ensure~\cite{ABV00} in particular the
existence of finitely many ergodic absolutely continuous (in
particular \emph{physical}) $f$-invariant probability
measures $\mu_1,\dots,\mu_k$ whose basins cover the manifold
Lebesgue almost everywhere.

We say that an $f$-invariant measure $\mu$ is an
\emph{equilibrium state} with respect to the potential $\log
J$, where $J=|\det Df|$, if $h_\mu(f)=\mu(\log J)$, that is
if \emph{$\mu$ satisfies the Entropy Formula}. Denote by
$\EE$ the family of all such equilibrium states.
It is not difficult to see that each physical measure in our
setting belongs to $\EE$.

We assume that \emph{ $\EE$ is formed by a unique
    absolutely continuous probability measure}.

\subsubsection{Conditions on the roof function}
\label{sec:condit-roof-functi}
We assume that $r:M\setminus\cS\to\RR^+$ has
\emph{logarithmic growth near $\cS$}: there exists
$K=K(\vfi)>0$ such that $r\cdot\chi_{B(\cS,\delta)}\le
K\cdot\big| \log d_{\delta}(x,\cS) \big|$ for all small
enough $\delta>0$, where $B(\cS,\delta)$ is the
$\delta$-neighborhood of $\cS$. We also assume that $r$ is
bounded from below by some $r_0>0$.

Now we can state the result on large deviations.

\begin{theorem}
  Let $X^t$ be a suspension semiflow over a non-uniformly
  expanding transformation $f$ on the base $M$, with roof
  function $r$, satisfying all the previouly stated conditions.

  Let $\psi:M_r\to\RR$ be continuous and
  $\nu=\mu\ltimes\leb^1$ be the induced invariant measure
  for the semiflow $X^t$, that is, for any $A\subset M_r$ we
  set $\nu(A)=\mu(r)^{-1}\int d\mu(x) \int_0^{r(x)}\! ds
  \,\chi_A(x,s).$ Let also $\lambda=\leb\ltimes\leb^1$ be
  the natural extension of volume to the space $M_r$.  Then
\begin{align*}
\limsup_{T\to\infty}\frac1T\log\lambda\left\{
z\in M_r: \left|
\frac1T\int_0^T \psi\left(X^t(z)\right)\, dt - \nu(\psi)
\right|>\epsilon 
\right\}<0.
\end{align*}
\end{theorem}

\subsubsection{Consequences for the geometric Lorenz flow}
\label{sec:conseq-geometr-loren}
Now consider a Lorenz geometric flow as constructed in
Section \ref{sec:geoLorenz} and let $F$ be the
one-dimensional map associated, obtained quotienting over
the leaves of the stable foliation, see Figure
\ref{fig:GeoLorenz}.  This map has all the properties stated
previously for the base transformation. The Poincar\'e
return time gives also a roof function with logarithmic
growth near the singularity line.

The exponentially slow recurrence property
\eqref{eq:expslowrec} depends on a delicate  combinatorial
argument for which the uniform expansion and the existence
of a unique singular point for the
one-dimensional map induced by the geometric Lorenz flow is
technically important. The proof of this property for
singular-hyperbolic attractors in general must deal with the
simultaneous presence of several singularities in the
one-dimensional map.

The uniform contraction along the stable leaves implies that
the \emph{time averages of two orbits on the same stable
  leaf under the first return map are uniformly close} for
all big enough iterates. 
If $P:S\to[-1,1]$ is the
projection along stable leaves
\begin{lemma}
  For $\vfi:U\supset\Lambda\to\RR$ continuous and bounded,
  $\epsilon>0$ and $\vfi(x)=\int_0^{r(x)}\psi(x,t)\,dt$,
  there exists $\zeta:[-1,1]\setminus\cS\to\RR$ with
  logarithmic growth near $\cS$ such that
$
\Big\{\big|\frac1n
S_n^{R}\vfi-\mu(\vfi)\big|>2\epsilon\Big\}$
is contained in 
$$
  P^{-1}\Big(
  \big\{\big|\frac1n S^f_n\zeta-\mu(\zeta)\big|>\epsilon\big\}
  \cup
  \big\{
  \frac1n S^f_n\big|\log\dist_\delta(y,\cS)\big|> \epsilon
  \big\}
\Big).
$$
\end{lemma}
Hence in this setting it is enough to study the quotient map
$f$ to get information about deviations for the Poincar\'e return map. 
Coupled with the main result we are then able to
deduce

\begin{corollary}
  Let $X^t$ be a flow on $\RR^3$ exhibiting a geometric
  Lorenz attractor with trapping region $U$. Denoting by
  $\leb$ the normalized restriction of the Lebesgue volume
  measure to $U$, $\psi:U\to\RR$ a bounded continuous
  function and $\mu$ the unique physical measure for the
  attractor, then for any given $\epsilon>0$
  \begin{align*}
      \limsup_{T\to\infty}\frac1T\log\leb\left\{
        z\in U: \left|
          \frac1T\int_0^T \psi\left(X^t(z)\right)\, dt - \mu(\psi)
        \right|>\epsilon 
      \right\}<0.
  \end{align*}
  Moreover for any compact $K\subset U$ such that $\mu(K)<1$
  we have
\begin{align*}
\limsup_{T\to+\infty}\frac1T \log \leb\Big( \left\{ x\in K :
  X^t(x)\in K, 0<t<T \right\} \Big) < 0.
\end{align*}
\end{corollary}

\subsubsection{Idea of the proof}
\label{sec:idea-proof}
We use properties of non-uniformly expanding
transformations, especially a large deviation bound recently
obtained~\cite{araujo-pacifico2006}, to deduce a large
deviation bound for the suspension semiflow reducing the
estimate of the volume of the deviation set to the volume of
a certain deviation set for the base transformation.

The initial step of the reduction is as follows.
For a continuous and bounded $\psi:M_r\to\RR$,  $T>0$
and $z=(x,s)$ with $x\in M$ and $0\le s < r(x) <\infty$,
 there exists the \textbf{lap number}
$n=n(x,s,T)\in\NN$ such that $S_{n}r(x)\le s+T < S_{n+1}
r(x)$,  and we can write
\begin{align*}
\int_0^T
\hspace{-0.2cm}\psi\big(X^t(z)\big)\,dt
=
\int_s^{r(x)}
\hspace{-0.6cm}\psi\big(X^t(x,0)\big)\,dt
&+\int_0^{T+s-S_{n}r(x)}
\hspace{-0.6cm}\psi\big(X^t(f^n(x),0)\big)\,dt
\\
&+\sum_{j=1}^{n-1} \int_0^{r(f^j(x))}
\hspace{-0.6cm}\psi\big(X^t(f^j(x),0)\big)\,dt.
\end{align*}
Setting $\vfi(x)=\int_0^{r(x)}\psi(x,0)\,dt$ we can rewrite
the last summation above as $S_n\vfi(x)$.  We get the
following expression for the time average
\begin{align*}
  \frac1T\int_0^T \hspace{-0.2cm}\psi\big(X^t(z)\big)\,dt
  =\frac1T S_n\vfi(x) 
  &
  -\frac1T\int_0^s\psi\big(X^t(x,0)\big)\,dt
  \\
  &+
  \frac1T\int_0^{T+s-S_n r(x)}
  \hspace{-0.6cm}\psi\big(X^t(f^n(x),0)\big)\,dt.
\end{align*}
Writing $I=I(x,s,T)$ for the sum of the last two integral
terms above, observe that for $\omega>0$, $0\le s<r(x)$ and
$n=n(x,s,T)$
\begin{align*}
  \left\{ (x,s)\in M_r : \left| \frac1T S_n\vfi(x) + I(x,s,T)
    -\frac{\mu(\vfi)}{\mu(r)}\right| > \omega\right\}
\end{align*}
is contained in
\begin{align*}
  \left\{ (x,s)\in M_r : \left| \frac1T S_n\vfi(x)
      -\frac{\mu(\vfi)}{\mu(r)}\right| > \frac\omega2 \right\}
\cup
\left\{(x,s)\in M_r: I(x,s,T) >\frac\omega2\right\}.
\end{align*}
The left hand side above is a \emph{deviation set for the
  observable $\vfi$ over the base transformation}, while the
right hand side will be \emph{bounded by the geometric
  conditions on $\cS$} and by a \emph{deviations bound for
  the observable $r$ over the base transformation}.

Analysing each set using the conditions on $f$ and $r$ and
noting that for $\mu$- and $\leb$-almost every $x\in M$ and every
$0\le s<r(x)$
\begin{align*}
 \frac{S_nr(x)}n\le
 \frac{T+s}{n}\le\frac{S_{n+1}r(x)}n
 \quad\text{so}\quad
 \frac{n(x,s,T)}T\xrightarrow[T\to\infty]{}\frac1{\mu(r)},
\end{align*}
we are able to obtain the asymptotic bound of the Main
Theorem.

Full details of the proof are presented
in~\cite{araujo2006a}.




\subsection{Decay of correlations for flows}
\label{sec:exponent-decay-corre-geomLorenzflow}


\subsubsection{Non-mixing flows and slow decay of correlations}
\label{sec:non-mixing-flows}

Let $f:M\to M$ be a diffeomorphism with an invariant
probability measure $\mu$ and consider the suspension flow
$X_f$ over $f$ with constant roof function $r\equiv1$. Then
the probability measure $\nu=\mu\times\leb$ on
$M\times[0,1)$ defines in a straightforward way a
$X_f$-invariant probability measure on $X_r$ \emph{which is
  NOT mixing}, whatever $f$ may be.

\begin{figure}[h]
  \psfrag{1}{$1$}\psfrag{0}{$0$}\psfrag{2}{$2$}
  \psfrag{3}{$3$}\psfrag{t}{$t$}
  \includegraphics[width=6cm]{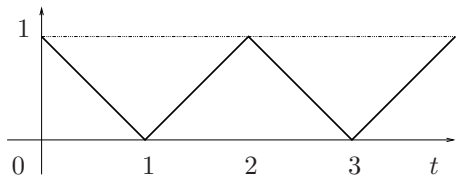}
  \caption{A correlation function for a non-mixing flow.}
  \label{fig:correl-functi-non}
\end{figure}

Indeed, consider $A=\pi\big(M\times[0,1/2)\big)$ and
$B=M_r\setminus A$ (recall that $\pi:M\times\RR\to X_r$ is
the projection defined in
Section~\ref{sec:large-deviations}).  Then the function
$t\mapsto\nu\big(A\cap X^{-t}B\big)$ for $t>0$ has the graph
as in Figure~\ref{fig:correl-functi-non} (here $X^{-t}$ is a
shorthand for $(X^t)^{-1}$, the inverse image of the map
$X^t$).

This system is clearly \emph{not}  mixing since the sawtooth
pattern in Figure~\ref{fig:correl-functi-non} goes on for
all positive $t$. Moreover this shows in particular that
this suspension flow is not even topologically mixing (see
below for the definition).

If however if $(X,f,\mu)$ is ergodic, then $\nu$ is
$X_f$-ergodic also: indeed, given $A\subset X_r$ such that
$(X_f^{t})^{-1}(A)=A$ for all $t>0$ (an $X_f$-invariant
set), then $A$ is saturated, i.e., $p\in A$ if, and only if,
$\cO_{X_f}(p)\subset A$; thus we may find $\hat A\subset X$
such that $A\cap\pi(X\times\{0\})=\pi(\hat A)$ is
$X_f^{1}$-invariant by construction (because $r\equiv1$),
$\hat A$ is $f$-invariant and $\nu(A)=\mu(\hat
A)\cdot\leb([0,1))$. Hence $\mu(\hat
A)\cdot\mu(X\setminus\hat A)=0$ by the ergodicity of
$(f,\mu)$ which implies that $\nu(A)\cdot\nu(X_r\setminus A)=0$.

In addition to the examples of non-mixing suspension flows,
which arguably can be characterized as very particular
cases, not all Axiom A mixing flows have exponential decay
of correlations: Ruelle \cite{ruelle1983} and Pollicott
\cite{pollicott84} exhibited suspensions semiflows with
piecewise constant ceiling functions over uniformly
expanding base dynamics, with arbitrarily slow decay rates
of correlations.

  The example from Ruelle is simple to describe: take the
  full shift on $2$ symbols $\sigma:\Sigma_2\to\Sigma_2$ and
  the roof function $r(\xi)=\lambda_0$ if $\xi_0=0$ and
  $r(\xi)=\lambda_1$ if $\xi_0=1$; where
  $\lambda_0,\lambda_1>0$ and $\lambda_0/\lambda_1$ is not
  rational. Take any equilibrium state $\mu$ for $\sigma$
  with respect to a H\"older continuous potential
  $\phi:\Sigma_2\to\RR$ and consider the induced probability
  $\nu=\mu\times\leb$ on $\{(\xi,s),0\le s<r(\xi)\}$. The
  suspension semiflow over $\sigma$ with roof function $r$
  does not have exponential decay of correlations for $\nu$.

Anosov \cite{An67} showed that geodesic flows for negatively
curved compact Riemannian manifolds are mixing and obtained
the \emph{Anosov alternative}: given a transitive volume
preserving Anosov flow, either it is mixing (with respect to
the volume measure), or it is a suspension of an Anosov
diffeomorphism by a constant roof function.  We note that
Bowen \cite{bowen76} showed that, if a mixing Anosov flow is
the suspension of an Anosov diffeomorphism, then it is
\emph{stably mixing}, that is, the mixing property remains
true for all nearby flows (which are Anosov also by the
structural stability of Axiom A flows).

Bowen also showed \cite{bowen76} that the class of $C^r$
Axiom A flows, $r\ge1$, admits a residual subset $\cR$ such
that for every $X\in\cR$ the spectral decomposition of
$\Omega(X)$ is formed by pairwise disjoint pieces
$\Omega_1\cup\dots\cup\Omega_k$ each of which is
\emph{topologically mixing}. That is, given any pair of open
sets $U, V$ in $\Omega_i$, there exists $T_0=T_0(U,V)>0$
such that $U\cap X^t(V)\neq\emptyset$ for all $t>T_0$.

\subsubsection{Exponential decay of correlations for
  hyperbolic flows}
\label{sec:recent-results}

Recently, a breakthrough was obtained by
Dolgopyat~\cite{Do98,dolgopyat98,dolgopyat2000}: smooth
($C^r$ with $r\ge7$) geodesic flows on manifolds of negative
curvature, under a non-integrability condition exhibit
exponential decay of correlations. Also Liverani
\cite{liverani2004} building on the work~\cite{Do98}
obtained exponential decay of correlations for $C^4$ contact
Anosov flows.

Using these ideas, applied to the particular case of a
suspension over uniformly expanding base dynamics, a
conjecture of Ruelle was proved by Pollicott
\cite{pollicott99}: on a mild (cohomological) condition on
the ceiling function, the decay of correlations for this
type of suspension flows is exponential for observables not
supported on the base.  This was extended by Baladi-Vall\'ee
\cite{BaVal2005}, clarifying the assumptions on the base and
on the ceiling function which suffice to obtain exponential
decay of correlations for suspension of one-dimensional
expanding maps. All these ideas were used, in a more
abstract setting, by Avila-Gouezel-Yoccoz \cite{AvGoYoc} to
obtain exponential decay of correlations for the
Teichm\"uller flow on flat surfaces.

Recently Field-Melbourne-T\"orok obtained~\cite{FMT} what
they call \emph{stability of rapid mixing} among Axiom A
flows, meaning that the correlation function
$C_t(\vfi,\psi)$ decays to zero faster than $t^{-k}$ for all
$k\in\NN$ when $t\to\infty$, for a $C^2$-open and
$C^r$-dense set of flows among the family of $C^r$ Axiom A
flows with $r\ge2$.

Luzzatto, Melbourne and Paccaut~\cite{LMP05} showed that the
physical measure for the geometric Lorenz flow, as presented
in Subsection~\ref{sec:geoLorenz}, is mixing.  The speed of
mixing for the Lorenz flow has been an open problem.

\subsubsection{Robust exponential decay of correlations
for a class of geometric Lorenz flows}
\label{sec:robust-exponent-deca}
The following result on exponential decay of
correlations for a suspension semiflow over a piecewise
expanding map with infinitely many branches is the
basis for a ongoing work \cite{ArVar} to find the rate
of decay of correlations for singular hyperbolic
attractors.

\begin{theorem}[Avila-Gouezel-Yoccoz \cite{AvGoYoc}]
  Let $Y_t$ be a good hyperbolic skew-product semi-flow on a
  space $\widehat{\Delta}_r$, preserving the probability
  measure $\bar\eta$. There exist constants $C>0$ and
  $\delta>0$ such that, for each pair of functions
  $\vfi,\psi\in C^1(\widehat{\Delta}_r)$, for all $t\geq 0$,
  \begin{align*}
  \left| \int \vfi\cdot \psi \circ Y_t \,d\bar\eta -\left( \int \vfi
  \,d\bar\eta\right) \left( \int \psi \,d\bar\eta \right) \right| \leq C
  \|\vfi\|_1 \|\psi\|_1e^{-\delta t}.
  \end{align*}
\end{theorem}

The meaning of "good" here will be explained below.

We show in \cite{ArVar} that an open class of geometric
Lorenz flows can be conjugated to semiflows in the above
setting, concluding robust exponential decay of correlations
for a wide class of singular flows.

\begin{theorem}[V.A.-P.Varandas \cite{ArVar}]
    \label{thm:Lorenz-semiconj}
  Given any compact $3$-manifold $M$, we
  can find an open subset $\U$ of $\X^3(M)$ such that each
  $X\in\U$ exhibits a geometric Lorenz flow which is
  smoothly semi-conjugated to a good hyperbolic skew-product
  semi-flow.
\end{theorem}

This result shows that the equilibrium point in the
geometric Lorenz flow actually \emph{helps} increase the
speed of decay of correlations, since the exponential decay
is in fact \emph{robust}. There are no examples of robust
exponential decay of correlations for Anosov flows! 

Now we explain what \textbf{\emph{good hyperbolic
    skew-product semiflow}} means, and how we relate a
geometric Lorenz flow to these semiflows.

\subsubsection*{Good hyperbolic skew-product semiflow}
\label{sec:good-hyperb-skew}

  We assume that $\cup_{\ell\in L} \De^{(\ell)}$ is an at
  most countable partition (Lebesgue modulo zero) of an open
  domain $\De$ of some manifold by open subsets and let
  $F:\cup_{\ell\in L} \De^{(\ell)} \to \De$ be a $C^r$
 \emph{uniformly expanding Markov map}, $r\ge2$, that is
\begin{enumerate}
\item $F:\De^{\ell}\to\De$ is a $C^{r}$
  diffeomorphism for every $\ell$;
\item there are $C>0$ and $0<\la<1$ such that 
  \begin{enumerate}
  \item for every inverse branch $h_n$ of $F^n$, with
    $n\ge1$,  $d(h_n(x),h_n(y)) \leq C \la^n d(x,y)$; and
  \item if $JF$ is the Jacobian of $F$ with respect to the
    Lebesgue measure, then $\log JF$ is a $C^1$ function and
    $\| D((\log JF)\circ h)\|_0\le C$ for every inverse
    branch $h$ of $F$.
  \end{enumerate}
\end{enumerate}
We denote by $\cH_n$ the family of inverse branches of
$F^n$. It is well known that $F$ admits an invariant
probability measure $\nu$ which is absolutely continuous
with respect to Lebesgue.

We say that the roof function $r$ is \emph{good} if
\begin{enumerate}
\item $r$ is bounded from below by some positive constant $r_0$;
\item there exists $C>0$ such that $\sup_{h\in\cH}\|D(r\circ
  h)\|_{0}\le C <\infty$;
\item it is not possible to write $r=v+u\circ F -u$ on
  $\De$, where $v: \De \to \R$ is constant on each
  $\De^{\ell}$ and $u:\Delta\to\RR$ is a $C^1$-function.
\end{enumerate}
The last cohomological condition corresponds to \textbf{uniform
non-integrability, or aperiodicity,} as defined by
Baladi-Vall\'ee adapted from the work of Dolgopyat.

We say that the roof function $r: \De \to \RR^+$ has
\emph{exponential tail} if there exists $\sigma_0>0$ such
that $\int e^{\sigma_0 r} d\nu <\infty$.

Let $F:\bigcup_{l} \Delta^{(l)} \to \Delta$ be a uniformly
expanding Markov map preserving a probability density
$\nu$. An \emph{hyperbolic skew-product} over $F$ is a map
$\widehat{F}$ from a dense open subset of an open domain
$\widehat{\Delta}$, to $\widehat{\Delta}$, satisfying
\begin{enumerate}
\item there exists a continuous map $\pi : \widehat{\Delta}
  \to \Delta$ such that $F\circ \pi = \pi \circ \widehat{F}$
  whenever both members of the equality are defined;
\item there is $\kappa>1$ such that, for all $w_1,w_2 \in
\widehat{\Delta}$ in the same leaf,
i.e. $\pi(w_1)=\pi(w_2)$, we have
$
  d(\widehat{F} w_1, \widehat{F} w_2) \leq \kappa^{-1}
  d(w_1,w_2).
$
\item there is a $\widehat F$-invariant probability measure
  $\eta$ on $\widehat{\Delta}$, giving full mass to
  $\widehat\Delta$;
\item there exists a smooth disintegration of $\eta$ along
  the stable leaves $\pi^{-1}(w), w\in\Delta$, as follows.
 
  There exists a family of probability measures
  $\{\eta_x\}_{x\in \Delta}$ on $\widehat{\Delta}$ which is
  a \textbf{disintegration} of $\eta$ over $\nu$:
  \begin{enumerate}
  \item $x\mapsto \eta_x$ is measurable;
  \item $\eta_x$ is supported on $\pi^{-1}(x)$, and
  \item for each measurable subset $A$ of
    $\widehat{\Delta}$ we have $\eta(A)=\int
    \eta_x(A)\,d\nu(x)$.
  \end{enumerate}
\end{enumerate}
Moreover, \emph{this disintegration is smooth}: we can
find a constant $C>0$ such that, for any open subset
$V\subset \bigcup \Delta^{(l)}$ and for each $u\in
C^1(\pi^{-1}(V))$, the function $\bar u : V \to \RR,
x\mapsto\bar u(x):=\int u(y)\,d\eta_x(y)$ belongs to
$C^1(V)$ and satisfies
  \begin{align*}
  \sup_{x\in V} \|D\bar u(x)\|
  \leq C \sup_{y\in \pi^{-1}(V)} \|Du(y)\|.
\end{align*}

Given $r : \cup_{\ell\in L} \De^{(\ell)} \to[r_0,+\infty)$
for some $r_0>0$ we define
\begin{align*}
  \widehat\De_r=\{(w,t): w \in \widehat\De, \; 0\leq t \le
  r(\pi(w))\}/\sim,
\end{align*}
where $\sim$ is the equivalence relation
$(w,r(\pi(w)))\sim(F(w),0)$.  Now we consider the suspension
semiflow $Y_t(w,s)=(w,s+t)$.

If $Y_t$ is a semiflow over a hyperbolic skew-product with a
good roof function which, moreover, has exponential tail,
then we say that $Y_t$ is a \emph{good hyperbolic
  skew-product semi-flow}.

We remark that, if $\eta$ is an $\widehat F$-invariant
probability measure so that $\int r\, d\eta<\infty$, then
$(Y_t)_t$ preserves the probability measure $ \bar\eta=(\eta
\otimes \leb)/ \int r \,d\eta.$

The following are the main steps of the proof of
Theorem~\ref{thm:Lorenz-semiconj}.

  \begin{itemize}
  \item obtain a robust $C^3$-smooth strong-stable foliation
    for the geometric Lorenz flow together with robust
    transitivity for the associated attractor;
  \item obtain a uniformly expanding $C^2$ Markov map $F$ as
    an induced map of the $C^2$ one-dimensional Lorenz
    transformation $f$;
  \item show that there is a related induced map $\hat F$
    from the Poincare return map $P$ to $S$ and a natural
    choice $r$ of roof function over $\hat F$ such that the
    original flow is conjugated to the semiflow over $\hat
    F$ with roof $r$;
  \item prove that $r$ has exponential tail, satisfies the
    aperiodicity condition and that the disintegration
    property holds for the right choice of measure on the
    semiflow over $\hat F$;
  \item check that each of the above steps are robust for
    $C^3$ close flows.
  \end{itemize}

  The crucial first item above depends on the
  inequality $\lambda_3>2\lambda_1+\lambda_2$ between
  the eigenvalues of the Lorenz-like singularity. This
  inequality does not hold in general and, in
  particular, does not hold for the flow of the Lorenz
  attractor.

\subsubsection{Rapid mixing for geometric Lorenz flows}
\label{sec:rapid-mixing-geometr}

It follows from~\cite{Melb09} that a $C^2$-open and
$C^\infty$-dense set of geometric Lorenz flows have
superpolynomial decay of correlations (in the sense of
\cite{dolgopyat98}); we also say that these flows are
\emph{rapidly mixing}.  It is likely, but unproven,
that this open and dense set includes the classical
Lorenz attractor.  

In a recently work \cite{AMV13} it was proved that
\emph{all $C^\infty$ geometric Lorenz attractors}
(including the concrete example of the attractor for
the Lorenz system of equations~\eqref{e-Lorenz-system})
satisfy this rapid mixing property.  More precisely,
let $\mathcal{U}$ denote the open set of $C^\infty$
vector fields having a geometric Lorenz attractor;
recall Section~\ref{3} for precise definitions.  Given
$X\in\mathcal{U}$, let $X^t$ denote the flow generated
by $X$ and let $\mu$ denote the unique physical measure
supported on the geometric Lorenz attractor.

\begin{theorem}{\cite[Theorem A]{AMV13}}\label{thm:rapid}
  Let $X\in\mathcal{U}$.  Then for all $\beta>0$, there
  exists $C>0$ and $k\ge1$ such that for all $C^k$
  observables $\vfi,\psi: \RR^3\to \RR$ and all $t>0$,
\[
\Big|
\int \vfi \; \psi\circ X^t \, d\mu -
\int \vfi\,d\mu  \int \psi\, d\mu\Big|
\le 
C\|\vfi\|_{C^k} \|\psi\|_{C^k}t^{-\beta}.
\]
\end{theorem}

This recent result encompasses all smooth geometric
Lorenz flows, including in particular the Lorenz
attractor given by (\ref{e-Lorenz-system}).

\subsection{Central limit theorem for the Lorenz flow}
\label{sec:central-limit-theore}

In \cite{HoMel} Holland and Melbourne, building on the work
\cite{MelNicol05} of Melbourne and Nicol, obtained the
Almost Sure Invariance Principle (ASIP) for \emph{Axiomatic
  geometrical Lorenz attractors} which, in turn, implies the
Central Limit Theorem.  More precisely, if $X^t$ is the
axiomatic geometric Lorenz flow with physical probability
measure $\mu$ and $\psi$ is a H\"older continuous function
(observable) on the manifold with zero mean
$\int\psi\,d\mu=0$, then there exists a Brownian motion
$W(t)$ with variance $\sigma^2>0$, and there is
$\epsilon>0$, such that
\begin{align*}
  \int_0^t\psi\circ X^s \, ds = W(t)+O(t^{\frac12-\epsilon})
  \quad\text{as}\quad t\to+\infty \quad\text{for
    $\mu$-almost all } x.
\end{align*}
This result, in turn, implies the Central Limit Theorem
(CLT): in
the same setting as above, for any interval $A\subset\RR$
\begin{align*}
  \mu\left\{
    x : \frac1{\sqrt{t}}\left(
       \int_0^t \psi\circ X^s \,ds - \mu(\psi)
      \right)\in A
    \right\}
    \xrightarrow[t\to+\infty]{}
    \frac1{\sigma\sqrt{2\pi}} \int_A e^{-\frac{s^2}{2\sigma}} \,ds;
\end{align*}
and the Law of the Iterated Logarithm
\begin{align*}
  \limsup_{t\to+\infty}\frac1{\sqrt{2t\log\log t}} \int_0^t
  \psi\circ X^s \,ds =\sigma \quad\text{$\mu$-almost everywhere}.
\end{align*}

In the recent work~\cite{AMV13} a stronger property was
obtained: the scalar ASIP holds for the time-$1$ map
$X^1$ of all smooth geometric Lorenz flows. This
information is used by the authors in~\cite{AMV13} to
prove the CLT for time-$1$ maps of geometric Lorenz
flows.

We note that, even for hyperbolic flows, the time-$1$
map is only partially hyperbolic. In general,
statistical results for the time-$1$ map are not
straightforward consequences of the continuous
versions.

\begin{theorem}{\cite[Theorem A]{AMV13}}\label{thm:CLT13}
  Let $X\in\mathcal U$.  Then there exists
  $k\ge1$ such that for all $C^k$ observables $\vfi:
  \RR^3 \to \RR$ there exists $\sigma\ge0$ such that
\begin{equation*}
  \frac{1}{\sqrt{n}} 
  \Bigg[\sum_{j=0}^{n-1} v\circ X^j \,
  -\,  n\int \vfi \;d\mu\Bigg] \xrightarrow{\mathcal D}
  \mathcal N(0,\sigma^2)
\end{equation*}
where the convergence is in distribution.
\end{theorem}

For the scalar ASIP we have, more precisely.

\begin{theorem}{\cite[Theorem A]{AMV13}}\label{thm:ASIP}
  Let $X\in\mathcal U$.  There exists $k\ge1$ such that for
  all $C^k$ observables $\vfi: \RR^3 \to \RR$ the ASIP holds for
  the time-$1$ map: passing to an enriched probability
  space, there exists a sequence $X_0,X_1,\ldots$ of iid
  normal random variables with mean zero and variance
  $\sigma^2$ (as in Theorem~\ref{thm:CLT13}), such that
\[
\sum_{j=0}^{n-1}\vfi\circ X^j=n\int \vfi\,d\mu+\sum_{j=0}^{n-1}
X_j+O(n^{1/4}(\log n)^{1/2}(\log\log n)^{1/4}),\; a.e.
\]
\end{theorem}

The ASIP implies the CLT and also the functional CLT
(weak invariance principle), and the law of the
iterated logarithm together with its functional
version, together with many other results;
see~\cite{PhilippStout75} for a comprehensive list.


\section{Some conjectures}
\label{sec:some-conject-open}

Singular-hyperbolic theory can be seen as an extension of
the theory of hyperbolicity, and so we may try to obtain the
same properties of hyperbolic sets in the setting of
singular-hyperbolic flows. Some of these results are already
known, for instance, the fact that every singular-hyperbolic
attractor is a homoclinic class, but their proof is,
usually, rather different from the usual ``hyperbolic
proof''. Other results are mostly wide open to research.
We present some of them here.

\subsection{Dimension theory, ergodic and statistical
  properties}
\label{sec:dimens-theory-ergodi}

Afraimovich and Pesin in~\cite{AP87} investigate the
dimensional properties of ``triangular maps'' which are a
class of maps generalizing the Poincar\'e first return map
$P$ of the geometric Lorenz model.

Concerning fractal dimensions of Lorenz attractors we
mention the results of Leonov \cite{Leo88,Leo01} together
with Bouichenko \cite{BoLeo89}.  The first contains explicit
formulas for the Lyapunov dimension of the Lorenz attractor
and, in the second, a simple upper bound on the Hausdorff
dimension of Lorenz attractors is given in terms of the
parameters of the Lorenz systems of
equations~\eqref{e-Lorenz-system}.  

\begin{conjecture}
  As in a hyperbolic attractor on a surface, the Hausdorff
  dimension of any singular-hyperbolic attractor on a
  $3$-manifold satisfies \emph{Bowen's formula}: it is the
  value $2+\gamma$ where $\gamma$ satisfies
  $P_{top}(\gamma\log|\det DX^1\mid E^s|)=0$, $P_{top}$ is
  the topological pressure of the attractor, and $E^s$ is
  the one-dimensional stable bundle over the attractor.
\end{conjecture}

In \cite{Young81} Young shows that the geometrical
Lorenz attractor can be approximated by horseshoes with
entropy close to that of the Lorenz attractor.

\begin{conjecture}
  It is possible to approximate the topological entropy of a
  singular-hyperbolic attractor by the topological entropy
  of horseshoes contained in the attractor.
\end{conjecture}

\section{Large deviations for the Lorenz flow}
\label{sec:large-deviat-geometr}

As explained in Section~\ref{sec:large-deviations}, for a
geometric Lorenz flow, the large deviations decay rate for
the volume/physical measure is exponential. That is, if we
set $\epsilon>0$ as an error margin and consider
\[
B_t=\Big\{
z: \Big|
\frac1t\int_{0}^{t}
  \psi\big( X^t(z) \big) 
-
\int\psi\,d\mu
\Big|>\epsilon
\Big\},
\]
then sufficient conditions were found, in terms of the base
transformation and the roof function, under which the
Lebesgue measure of $B_t$ decays to zero exponentially fast,
i.e., whether there are constants $C,\xi>0$ such that
\begin{align*}
\leb\big( B_t\big) \le C e^{-\xi t}
\quad\mbox{for all}\quad t>0.
\end{align*}
We observe that in this setting Lebesgue measure or volume
is \emph{not an invariant measure}.

In \cite{melnicol07} Melbourne and Nicol and in
\cite{BellYoung08} Rey-Bellet and Young obtained large
deviations principles for \emph{invariant measures} in the
same setting, including subexponential or polynomial bounds
on large deviations depending on the properties of the base
transformation.

\begin{conjecture}
  These results are also true for general
  singular-hyperbolic attractors and should be true, under
  some mild conditions, for singular-hyperbolic attracting
  sets as well.
\end{conjecture}

\subsection{Decay of correlations}
\label{sec:decay-correl}

As explained in Section~\ref{sec:Hypflows-decay}, it
has recently been obtained an example of a geometric
Lorenz flow with robust exponential decay of
correlations for all flows sufficiently $C^2$ close.

A main open question is to find the rate of decay of
correlation for the original Lorenz flow, and more in
general for singular hyperbolic attractors.  In the
direction of finding an appropriate functional analytic
framework to face the problem, some avances was made
recently in \cite{Butterley12}.

 Some other natural questions are as follows.

\begin{conjecture}
  Non-hyperbolic robustly {mixing} flows in
  three-dimensional manifolds have robust exponential decay
  of correlations.
\end{conjecture}

The rate of decay should depend continuously on the system considered.

\begin{conjecture}
  Fixing the observable $\psi$, the rate of decay of
  correlations should depend continuously on the flow $X^t$
  in a neighborhood of a singular-hyperbolic attractor.
\end{conjecture}

Moreover, since the argument leading to robust
exponential decay crucially depends on the smoothness
of the stable foliation, which cannot be obtained
robustly for globally hyperbolic flows, we conjecture
the following.

\begin{conjecture}
  There are no Anosov flows with robust exponential decay of
  correlations.
\end{conjecture}

However, we should be able to obtain such smooth
foliations for flows whose limit set is hyperbolic.

\begin{conjecture}
  There are open sets of $C^2$ Axiom A flows on compact
  manifolds exhibiting robust exponential decay of
  correlations.
\end{conjecture}

For the definition of Axiom A the reader should consult
a standard reference on hyperbolic dynamics, e.g. \cite{Sm67}.

\subsection{Central Limit Theorem}
\label{sec:central-limit-theore-1}

We note that the ``time-$1$'' map $X_1$ of any flow $X_t$ on
a hyperbolic or singular-hyperbolic attractor is a partially
hyperbolic diffeomorphism. In general, limit theorems for
diffeomorphisms given as time-$t$ maps of flows are harder
to obtain.   A very general result was obtained by
Melbourne and T\"orok~\cite{MelTor02} under some assumptions
on the decay of correlations for the flow.   These
ideas can be adapted to prove that the strong mixing
properties for the $C^2$-open subset of geometric Lorenz
attractors imply (robust) limit theorems for the
corresponding time-one maps. More precisely we pose the
following:

\begin{conjecture}\label{conj:CLT-Time1}
  Let $\cU\subset \fX^s(M)$ be the open family of vector
  fields for which exponential decay of correlations is
  verified, and denote by $(X_t)_t$ the flow generated by
  $X\in \cU$. For all but countably many values of
  $t\in\mathbb R$ the time-$t$ map $X_t$ the following
  Central Limit Theorem holds: for any $\varphi: \De_r \to
  \RR$ in $L^\infty(\De_r)$ there exists
  $\sigma=\sigma(\varphi)>0$ such that
\begin{equation*}\label{CLT1}
  \frac{1}{\sigma \sqrt{n}} 
  \Bigg[\sum_{j=0}^{n-1} \varphi(X_{tn}) - 
  \int \varphi \;d\mu\Bigg] \xrightarrow{\mathcal D}
  \mathcal N(0,1).
\end{equation*}
\end{conjecture}

\subsection{Thermodynamical formalism}
\label{sec:thermodyn-formal}

The thermodynamical formalism was first developed for
(uniformly) hyperbolic diffeomorphisms, borrowed from
statistical mechanics by Bowen, Ruelle and Sinai (among
others, see e.g.
\cite{Bo75,BR75,Ru89,ruelle2004,ellis06,BDV2004}).
This was extended to hyperbolic flows by Bowen and Ruelle in
\cite{BR75}. The classical theory relies heavily  on the
coding of basic pieces of hyperbolic dynamics by subshifts
of finite type, for which many tools are available to study
in fine detail the relations among its invariant measures.
Recently most of this theory was extended to countable
shifts by Gurevich \cite{GurevSavch98}, Sarig
\cite{Sarig99,Sarig2006} and many others.

The extension of this theory for singular-hyperbolic
attractors faces several difficulties: these attractors are
modelled by a suspension semiflow whose base transformation
is a H\"older-$C^1$ piecewise expanding but
\emph{non-Markov} map, and the roof function is
unbounded. In the hyperbolic flow case, the corresponding
suspension semiflow has a piecewise expanding Markov map as
the base transformation and the roof function is continuous
and bounded. In the singular-hyperbolic case,
neither the thermodynamical formalism is complete
for the base transformation, nor  is it clear how to proceed
with unbounded roof functions, which imply an extra
restriction of integrability on the observables with respect
to invariant measure for the base transformation.

Hope of solving this problem in the near future is provided
by recent advances in the construction of a thermodynamical
formalism for non-uniformly expanding transformations by
Oliveira, Viana, Senti, Pesin, Varandas, Bruin, Todd,
Pinheiro
\cite{OliVi2005,pesin2006emm,VarVia09,BrTdd07,BrTdd08,Pinho2011}.

Recently, Leplaideur and Pinheiro in \cite{LepPinh13}
obtain the first results of a thermodynamic formalism
for a two-dimensional map representing the Poincar\'e
first return map of an expanding geometric Lorenz
attractor. Namely, they prove the existence of unique
equilibrium state for any H\"older continuous potential
on the attracting set of this two-dimensional map.

\begin{conjecture}
  The same result on existence and uniqueness of
  equilibrium states can be extended to all
  singular-hyperbolic attractors in three-dimensional
  manifolds.
\end{conjecture}

Ongoing work \cite{GT09,PT09} on the study of
equilibrium states for multiples of the logarithm of
the derivative, for suspension flows over
transformations resembling the Lorenz one-dimensional
transformations, will enable the extension of this
results from discrete dynamics to the flow of
Lorenz-like attractors.

\begin{conjecture}
  It is possible to build a thermodynamical formalism
  for Rovella-like and singular-hyperbolic attractors
\end{conjecture}

\subsection{Higher dimensional singular flows}
\label{sec:higher-dimens-singul}

An example of a higher-dimensional invariant robust
attractor with multidimensional expanding directions was
given by Bonatti, Pumarino and Viana in~\cite{BPV97}, which
we present below.

\subsubsection{Singular-attractor with arbitrary number of
  expanding directions}
\label{sec:singul-attract-with}

Consider a ``solenoid'' constructed over a uniformly expanding map
$f:\TT^k\to\TT^k$ of the $k$-dimensional torus, for some $k\ge2$.
That is, let $\DD$ be the unit disk on $\RR^2$ and consider a smooth
embedding $F:\TT^k\times\DD\to\TT^k\times\DD$ of $N=\TT^k\times\DD$
into itself, which preserves and contracts the foliation
$$\F^s=\big\{\{z\}\times\DD: z\in\TT^k\big\},$$ and moreover the natural
projection $\pi:N\to\TT^k$ on the first factor conjugates $F$ to
$f$: $\pi\circ F=f\circ\pi$.

\begin{figure}[htbp]
\includegraphics[width=7.5cm]{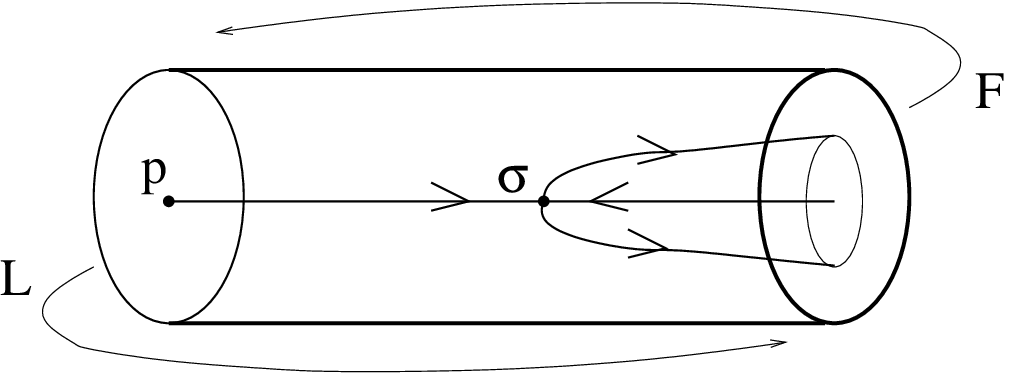}
\caption{A sketch of the construction of a
  robust singular-attractor in higher dimensions}
\label{fig:sing-attractor-4d}
\end{figure}

Now consider the linear flow over $M=N\times[0,1]/\sim$
given by the vector field $X=(0,1)$ on $TN\times\RR$ where
we make the identification $(x,0)\sim(x,1)$ for all $x\in
N$. Modify the flow on a cylinder $U\times\DD\times[0,1]$
around the orbit of a point $p=(z,0)\in N$, where $U$ is a
neighborhood of $z$ in $\TT^k$, in such a way as to create a
hyperbolic singularity $\sigma$ of saddle-type with
$k$-expanding and $3$ contracting eigenvalues, as depicted
in Figure~\ref{fig:sing-attractor-4d}.

This modified flow defines a transition map $L$ from
$\Sigma_0=\TT^k\times\{0\}$ to $\Sigma_1=\TT^k\times\{1\}$
which through the identification
given by $(w,1)\sim_F (F(w),0)$ defines the return map to
the global cross-section $\Sigma_0$ of a flow $Y$ on the
space $M^F=M/\sim_F$.

In~\cite{BPV97} it is shown that, if the expanding rate of $f$ is
sufficiently big, then the set
\[
\Lambda= \bigcup_{T>0} \overline{ \bigcap_{t>T} Y_t(\Sigma_0) }
\]
is a robust partially hyperbolic attractor with
singularities.

\subsubsection{The notion of sectionally hyperbolic sets}
\label{sec:section-expand-attra}

Metzger and Morales in \cite{MeMor06} introduced the
notion of sectionally expanding
or\emph{sectional-hyperbolic set} in a manifold of
arbitrary finite dimension. This notion encompasses
that of singular-hyperbolic sets in $3$-manifolds as a
particular case.

We say that a compact invariant set $\Lambda$ for a
flow, generated by a vector field $X\in\fX^1(M)$ on a
compact finite dimensional manifold $M$, is
\emph{sectional-hyperbolic} if it is partially
hyperbolic and the central direction expands uniformly
the area along any two-dimensional subspace. More
precisely, the tangent bundle over $\Lambda$ admits a
$DX^t$-invariant and dominated splitting $T_\Lambda M=
E^s\oplus E^c$, such that there are $C,\lambda>0$
satisfying for every $x\in\Lambda$ and $t>0$
\begin{itemize}
\item $E^s$ is uniformly contracted: $\|DX^t\mid E^s_x\|\le
  C e^{-\lambda t}$;
\item $E^c$ is sectionally expanded: for every
  bidimensional subspace $F_x$ contained in $E^c_x$ we have
  $|\det (DX^t\mid F_x)| \ge C e^{\lambda t}$.
\end{itemize}

Similarly to the notion of singular-hyperbolicity, robust
attractors in higher dimensional manifolds need not be
sectional-hyperbolic, as the example of Turaev and Shil'nikov
in \cite{ST98} shows.

The results in Section~\ref{sec:absence-sinks-source} have
a precise counterpart for sectional-hyperbolic attractors
with very similar proofs. It is then natural to try to
extend the three-dimensional results to this more general
setting.

\begin{conjecture}
  All the results obtained for singular-hyperbolic
  attractors should hold true for sectional-hyperbolic
  attractors in any dimension.
\end{conjecture}

\section{Appendix: chaotic systems}
\label{sec:chaotic-systems}

We distinguish between forward and backward sensitive
dependence on initial conditions. We say that an invariant
subset $\Lambda$ for a flow $X^t$ is \emph{future chaotic
  with constant $r>0$} if, for every $x\in\Lambda$ and each
neighborhood $U$ of $x$ in the ambient manifold, there
exists $y\in U$ and $t>0$ such that
$\dist\big(X^t(y),X^t(x)\big)\ge r$. Analogously we say that
$\Lambda$ is \emph{past chaotic with constant $r$} if
$\Lambda$ is future chaotic with constant $r$ for the flow
generated by $-X$. If we have such \emph{sensitive
  dependence both for the past and for the future}, we say
that $\Lambda$ is \emph{chaotic}. Note that in this language
sensitive dependence on initial conditions is weaker than
chaotic, future chaotic or past chaotic conditions.

An easy consequence of chaotic behavior is that it prevents
the existence of sources or sinks, either attracting or
repelling equilibria or periodic orbits, inside the
invariant set $\Lambda$. Indeed, if $\Lambda$ is future
chaotic (for some constant $r>0$) then, were it to contain
some attracting periodic orbit or equilibrium, any point of
such orbit (or equilibrium) would admit no point in a
neighborhood whose orbit would move away in the
future. Likewise, reversing the time direction, a past
chaotic invariant set cannot contain repelling periodic
orbits or repelling equilibria. As an almost reciprocal we
have the following.

\begin{lemma}
  \label{le:nonchaotic-interio}
  If $\Lambda=\cap_{t\in\RR}X^t(U)$ is a compact isolated
  proper subset for $X\in\fX^1(M)$ with isolating
  neighborhood $U$ and $\Lambda$ is \emph{not future
    chaotic} (respective not past chaotic), then
  $\Lambda^-_X(U):=\overline{\cap_{t>0}X^{-t}(U)}$
  (respective $\Lambda^+_X(U):=\overline{\cap_{t>0}X^t(U)}$)
  has non-empty interior.
\end{lemma}

\begin{proof}
  If $\Lambda$ is not future chaotic, then for every $r>0$
  there exists some point $x\in\Lambda$ and a neighborhood
  $V$ of $x$ such that $\dist\big( X^t(y), X^t(x)\big)<r$
  for all $t>0$ and each $y\in V$. If we choose
  $0<r<\dist(M\setminus U,\Lambda)$ (we note that if
  $\Lambda=U$ then $\Lambda$ would be open and closed, and
  so, by connectedness of $M$, $\Lambda$ would not be a
  proper subset), then we deduce that $X^t(y)\in U$, that
  is, $y\in X^{-t}(U)$ for all $t>0$, hence $V\subset
  \Lambda^-_X(U)$. Analogously if $\Lambda$ is not past
  chaotic, just by reversing the time direction.  \qed
\end{proof}

In particular \emph{if an invariant and isolated set $\Lambda$
  with isolating neighborhood $U$ is given such that the volume of
  both $\Lambda_X^+(U)$ and $\Lambda^-_X(U)$ is zero, then
  $\Lambda$ is chaotic}.

Sensitive dependence on initial conditions is part of many
definitions of \emph{chaotic behavior} in the literature, see
e.g.~\cite{devaney1989}. It is an interesting fact that
sensitive dependence is a consequence of another two common
features of most systems considered to be chaotic: existence
of a dense orbit and existence of a dense subset of periodic
orbits. 

\begin{proposition}
  \label{pr:transdensperiodic-sensitive}
  A compact invariant subset $\Lambda$ for a flow $X^t$ with
  a dense subset of periodic orbits and a dense (regular and
  non-periodic) orbit is chaotic, in the sense defined above.
\end{proposition}

A short proof of this proposition can be found
in~\cite{BBCDS92}. An extensive discussion of this and
related topics can be found in~\cite{GlasWeiss93}.

\subsection{Robust chaoticity, volume hyperbolicity and
  physical measure}
\label{sec:robust-chaotic-parti}

Here we show that robust chaotic behavior of an attractor,
under mild conditions, is equivalent to
singular-hyperbolicity and ensures the existence of a
physical measure. This is proved by noting a series of
consequences of the previous results

Here, by ``robust chaotic attractor'' $\Lambda=\Lambda_X(U)$
in a trapping region $U$ of a vector field $X$ we mean that,
for all close enough vector fields $Y$ to $X$ in the $C^1$
topology, the corresponding maximal invariant subset
$\Lambda_Y(U)$ is also chaotic.

Going through the proof of
Theorem~\ref{thm:robust-attractor-sing-hyp} in~\cite{MPP04}
we can see that the arguments can be carried through
assuming that
\begin{enumerate}
\item $\Lambda$ is an attractor for $X$ with isolating
  neighborhood $U$ such that every equilibria in $U$ is
  hyperbolic with no resonances;
\item there exists a $C^1$ neighborhood $\cU$ of $X$ such
  that for all $Y\in\cU$ every periodic orbit and
  equilibria in $U$ is hyperbolic of saddle-type.
\end{enumerate}
The condition on the equilibria amounts to restricting
the possible three-dimensional vector fields in the above
statement to an open a dense subset of all $C^1$ vector
fields. Indeed, the hyperbolic and no-resonance condition on
a equilibrium $\sigma$ means that:
\begin{itemize}
\item either $\lambda\neq\Re(\omega)$ if the eigenvalues of
  $DX(\sigma)$ are $\lambda\in\RR$ and
  $\omega,\overline\omega\in\CC$;
\item or $\sigma$ has only real eigenvalues with different
  norms.
\end{itemize}
Indeed, conditions (1) and (2) ensure that no bifurcations
of periodic orbits or equilibria leading to sinks or sources
are allowed for any nearby flow in $U$. This implies, by now
standard arguments, that the flow on $\Lambda$ must have a
dominated splitting which is \emph{volume hyperbolic}: both
subbundles of the splitting must contract/expand volume; see
e.g. \cite{AraPac2010}. For a $3$-dimensional flow one of
the subbundles is one-dimensional, and so we deduce
singular-hyperbolicity either for $X$ or for $-X$.  If
$\Lambda$ has no equilibria, then $\Lambda$ is uniformly
hyperbolic. Otherwise, it follows from the arguments
in~\cite{MPP04} that all singularities of $\Lambda$ are
Lorenz-like and this shows that $\Lambda$ must be
singular-hyperbolic for $X$.

We note that condition (2) above is a consequence of any one
of the following assumptions a neighborhood of $X$ and on
the neighborhood $U$ in $M$:
\begin{description}
\item[\textbf{robust chaoticity}] for every $Y\in\cU$ the maximal
  invariant subset $\Lambda_Y(U)$ is chaotic;
\item[\textbf{zero volume and future chaoticity}] for every
  $Y\in\cU$ the maximal invariant subset $\Lambda_Y(U)$ has
  zero volume and is future chaotic;
\item[\textbf{zero volume and robust positive Lyapunov
  exponent}] for every $Y\in\cU$ the maximal
  invariant subset $\Lambda_Y(U)$ has zero volume and
  there exists a full Lebesgue measure subset $P_Y$ of $U$
  such that
  \begin{align}
    \label{eq:positiveLyap}
    \limsup\frac1n\sum_{i=0}^{n-1}\log\|DY^i_x\|>0, \quad
    x\in P_Y.
  \end{align}
\end{description}
The following result of  Ma\~n\'e analogous to
Theorem~\ref{thm:robust-attractor-sing-hyp} in~\cite{Man82}
also follows from the absence of sinks and sources for all
$C^1$ close diffeomorphisms in a neighborhood of the
attractor.
\begin{theorem}
  \label{thm:robust2d-hyp}
  Robust attractors for surface diffeomorphisms are hyperbolic.
\end{theorem}

Extensions of these results to higher dimensions for
diffeomorphisms, by Bonatti, Diaz and Pujals in~\cite{BDP},
show that robust transitive sets always admit a volume
hyperbolic splitting of the tangent bundle.  Vivier in
\cite{Viv03} extends previous results of Doering~\cite{Do87}
for flows, showing that a $C^1$ robustly transitive vector
field on a compact boundaryless $n$-manifold, with $n\ge3$,
admits a global dominated splitting. Metzger and Morales
extend the arguments in \cite{MPP04} to homogeneous vector
fields (inducing flows allowing no bifurcation of critical
elements, i.e. no modification of the index of periodic
orbits or equilibria) in higher dimensions leading to the
concept of $2$-sectional expanding attractor
in~\cite{MeMor06}.

The preceding observations allows us to deduce that robust
chaoticity is a sufficient condition for
singular-hyperbolicity of a generic attractor; see
\cite{AraPac2010}.

\begin{corollary}
  \label{cor:robust-chaotic-sing-hyp}
  Let $\Lambda$ be an attractor for $X\in\fX^1(M^3)$ such
  that every equilibrium in its trapping region is
  hyperbolic with no resonances. Then $\Lambda$ is
  singular-hyperbolic if, and only if, $\Lambda$ is robustly
  chaotic.
\end{corollary}

This means that \emph{if we can show that arbitrarily close
  orbits, in an isolating neighborhood of an attractor, are
  driven apart, for the future as well as for the past, by
  the evolution of the system, and this behavior persists
  for all $C^1$ nearby vector fields, then the attractor is
  singular-hyperbolic}.

To prove the necessary condition on
Corollary~\ref{cor:robust-chaotic-sing-hyp} we use the
concept of expansiveness for flows, and through it show that
singular-hyperbolic attractors for $3$-flows are robustly
expansive and, as a consequence, robustly chaotic
also. This is \cite[Theorem A]{APPV} whose
  proof is presented also in \cite[Chapter 7, Section
  7.2]{AraPac2010}.

We recall the following conjecture of Viana, presented in~\cite{Vi98}
\begin{conjecture}\label{conj:viana}
  If an attracting set $\Lambda(U)$ of smooth map/flow has a
  non-zero Lyapunov exponent at Lebesgue almost every point
  of its isolated neighborhood $U$ (i.e. it
  satisfies~\eqref{eq:positiveLyap} with $P_Y\subset U$),
  then it admits some physical measure.
\end{conjecture}
From the preceding results and observations we can give a
partial answer to this conjecture for $3$-flows in the
following form.

\begin{corollary}
\label{cor:weakconjecture}
Let $\Lambda_X(U)$ be an attractor for a flow $X\in\fX^1(M)$
such that
\begin{itemize}
\item the divergence of $X$ is negative in $U$;
\item the equilibria in $U$ are hyperbolic with no resonances;
\item there exists a neighborhood $\cU$ of $X$ in $\fX^1(M)$
  such that for $Y\in\cU$ one has \eqref{eq:positiveLyap}
  almost everywhere in $U$.
\end{itemize}
Then there exists a neighborhood $\cV\subset\cU$ of $X$ in
$\fX^1(M)$ and a dense subset $\cD\subset\cV$ such that
\begin{enumerate}
\item  $\Lambda_Y(U)$ is singular-hyperbolic for all $Y\in\cV$;
\item there exists a physical measure $\mu_Y$ supported in
  $\Lambda_Y(U)$ for all $Y\in\cD$.
\end{enumerate}
\end{corollary}
Indeed, item (2) above is a consequence of item (1), the
denseness of $\fX^2(M)$ in $\fX^1(M)$ in the $C^1$ topology,
together with Theorem~\ref{srb} and the observation
following its statement.

Item (1) above is a consequence of
Corollary~\ref{cor:robust-chaotic-sing-hyp} and the
observations of Section~\ref{sec:absence-sinks-source},
noting that negative divergence on the isolating
neighborhood $U$ ensures that the volume of $\Lambda_Y(U)$
is zero for $Y$ in a $C^1$ neighborhood $\cV$ of $X$.

\def\cprime{$'$}

 \bibliographystyle{ws-ijbc}

\end{document}